\documentclass[a4paper, 11pt]{amsart}
\usepackage{geometry}
\usepackage{amsmath,amssymb,amsfonts,amsthm} % need to fix margins if changing font siz
\usepackage{color}
\usepackage{fullpage}
\usepackage[all]{xy}
\usepackage[small,nohug,heads=vee]{diagrams}
\diagramstyle[labelstyle=\scriptstyle]

\DeclareMathOperator{\Supp}{Supp}
\newtheorem{thm}{Theorem}
  \newtheorem{cor}[thm]{Corollary}
\newtheorem{lem}[thm]{Lemma}
\newtheorem{prop}[thm]{Proposition}
\newtheorem{conjec}[thm]{Conjecture}
\theoremstyle{definition}
\newtheorem{prob}[thm]{Problem}
\newtheorem{defn}[thm]{Definition}
\theoremstyle{remark}
\newtheorem{rem}[thm]{Remark}
\newtheorem{claim}[thm]{Claim}
\theoremstyle{remark}

\numberwithin{equation}{section}
\numberwithin{thm}{section}
\newcommand{\HK}{hyper-K\"ahler}
\newcommand{\PP}{{\mathbb{P}}}

\newcommand{\ZZ}{{\mathbb{Z}}}
\newcommand{\oo}{{\mathcal O}}

\DeclareMathOperator{\supp}{supp}

\DeclareMathOperator{\mult}{mult}
\title{On irreducible symplectic 4-folds numerically
          equivalent to $(K3)^{[2]}$ }
\date{\today}
\author[G.~Kapustka]{Grzegorz Kapustka}
\address{University of Z\"urich, Jagiellonian University Cracow}
\email{grzegorz.kapustka@uj.edu.pl}
\begin{document}

\begin{abstract} 
We address the problem of classification of hyper-K\"ahler fourfolds with $b_2=23$.
In particular we prove some special cases of the Conjecture of O'Grady about \HK{} $4$-folds numerically
 equivalent to the Hilbert scheme of two points on a $K3$ surface.
\end{abstract}

\thanks{
Partially
supported by MNiSW grant N N201 388834. The author is supported by
the Foundation for Polish Science (FNP)}
\maketitle
\section{Introduction}
In this work we study the geometry of \HK{} fourfolds. 
\begin{defn}  A complex manifold $X$  of even dimension is called irreducible holomorphic symplectic or hyper-K\"ahler if 
\begin{enumerate}
\item $X$ is compact and K\"ahler, 
\item $X$ is simply connected, and 
\item $H^0(X,\Omega^2_X)$ is spanned by an everywhere non-degenerate two-form $\rho$.
\end{enumerate}
\end{defn}
Recall that hyper-K\"ahler manifolds are among the building blocks of K\"ahler fourfolds with trivial first Chern class. 
\begin{thm}(\cite{B},\cite{Bogomolov}) Let $X$ be a compact K\"ahler manifold with trivial Ricci curvature.
There exists an \'{e}tale finite cover $X'$ of $X$ isomorphic as a K\"ahler manifold to the finite product
$$T\times \Pi V_i \times \Pi X_i$$ 
where $T$ is a complex torus, $V_i$ are compact K\"ahler varieties that are simply connected with holonomy group $SU(m_i)$, with $\dim V_i=m_i$, and $X_i$ are compact K\"ahler varieties that are simply connected with holonomy group $Sp(m_i)$ and $\dim X_i=2m_i$.
\end{thm}
We call $V_i$ Calabi--Yau varieties and we are interested in the hyper-K\"ahler manifolds $X_i$. 
In dimension $2$ hyper-K\"ahler manifolds are called $K3$ surfaces.
It is known that all $K3$ surfaces are deformation equivalent.
In higher dimension much less is known about hyper-K\"ahler manifolds and there are no general classification results.
In fact we know in each even dimension $2n$ two families: the deformation of the Hilbert scheme  $S^{[n]}$ of $n$ points on a K3 surface $S$ called of \emph{type $K3^{[n]}$} and deformations of the Hilbert scheme $K_n(T)$ of $n+1$ points summing to $0$ on an abelian surface $T$. 
Moreover, two sporadic families of such varieties constructed by O'Grady, one in dimension $6$
 and the other in dimension $10$.
 % (see \cite{Ogrady1}).

%\begin{rem}
%More compact connected $4n$-dimensional Riemannian manifold $(M,g)$ is called hyperk\"ahler (irreducible %hyperk\"ahler) if its holonomy is contained in (equals) $Sp(n)$.
%Irreducible symplectic manifolds with a fixed K\"ahler class and compact irreducible hyperk\"ahler manifolds are incarnation of the same object.
%\end{rem} 

A fundamental fact about hyper-K\"ahler manifolds (proved in \cite{B}, \cite{Fujiki}) is that the intersection form on 
$H^2(X,\ZZ)$ for a
hyper-K\"ahler manifold $X$ is defined by the Beauville-Bogomolov (B-B) quadratic form $q \colon H^2(X, \ZZ) \to \ZZ$ that is primitive, of 
signature $(3, b_2 -3)$. We define it (see \cite{B}) through the Fujiki relation $$\int \alpha^{2r}= f_X q(\alpha)^r$$ for each $\alpha\in H^2(X,\ZZ)$, where $\dim(X) = 2r$.
We call $f_X>0$ the Fujiki constant. It is an open problem to bound the possible values of $f_X$ and find the possible quadric forms $q$ that can occur for 
hyper-K\"ahler manifolds.

There are known restrictions on the number of families of hyper-K\"ahler manifolds.
In dimension $4$ Guan proved that if $X$ is a hyper-K\"ahler fourfold then either the second Betti number is bounded $3\leq b_2 \leq 8$, or $b_2 = 23$. 
Note that the manifolds $S^{[2]}$ have the second Betti number $23$ and the manifolds $K_2(T)$ have $7$. 
%Huybrechts \cite{Hu} observes that there are finitely many deformation types of hyper-K\"ahler manifolds X of dimension 2r such that there exists $\alpha \in H^2(X, \ZZ)$ with $q(\alpha)>0$ and
%$\alpha^{2r}$ bounded. As a corollary, given a real number there are finitely many families of hyper-K\"ahler manifolds with
%$f_X \leq M$, $$min\{ q(\alpha)|\alpha \in H^2(X,\ZZ), q(\alpha)>0 \}\leq M .$$
%More restriction was found \cite{Kamenova}
It is expected that there are only a finite number of families of hyper-K\"ahler manifolds in each dimension.
See \cite{Kamenova} for a survey about known results concerning the finiteness of hyper-K\"ahler manifolds.

In this paper we address the problem of classification of \HK{} fourfolds with $b_2=23$. It is expected that all such fourfolds are deformations of the Hilbert scheme of two point on a $K3$ surface (i.e.~ are of type $K3^{[2]}$).
This problem is studied by O'Grady \cite{O} in a more special case. He aim to classify \HK{} fourfolds with the same Fujiki invariant and the same B-B form that the Hilbert scheme of two points on a $K3$ surface. We call such fourfolds \emph{numerical $K3^{[2]}$}, note that the Betti number $b_2=23$ in this cases. 
%Let us be more precise:we call two \HK{} manifolds $M_1$, $M_2$ of dimension
%$2n$ are numerically equivalent if there exists an isomorphism of
%abelian groups $$\psi\colon H^2(M_1, \mathbb{Z})\rightarrow
%H^2(M_2,\mathbb{Z})$$ such that $\int_{M_1} \alpha^{2n} = \int_{M_2}
%\psi(\alpha)^{2n}$ for all $\alpha \in H^2(M_1, \mathbb{Z})$. 
In \cite{O} O'Grady stated the following:
\begin{conjec}\label{conj} Hyper-K\"ahler fourfolds that are numerical $K3^{[2]}$ are of type $K3^{[2]}$.
\end{conjec}
In \cite{O} O'Grady presents a program to aproach Conjecure \ref{conj}.
The idea is to consider the behavior of maps to $\PP^5$ defined by an ample divisors $H$ with minimal $q(H)=2$ on such manifolds.
The image of such a map is a hypersurface in $\PP^5$ and the program consist in finding which hypersurfaces of degree $6\leq d\leq 12$ in $\PP^5$ can be birational to a
\HK{} fourfold (see Section \ref{OG}).
The aim of this paper is contribute to this program by doing some progress in the O'Grady program concerning this Conjecture \ref{conj}. 
 In particular we describe projective models of \HK{} fourfolds in $\PP^5$ that are the possible counterexamples to Conjecture \ref{conj}. 
 %Note that in \cite{K} we show that if $d=12$ the projective models are strongly related with special sextic hypersurfaces in $\PP^5$ called EPW sextics.

\section{The O'Grady program}\label{OG}
The fact that all K3 surfaces are deformation equivalent was first proved by Kodaira. 
In \cite{LP} Le Potier gave another proof of this fact that follows the following steps:
\begin{enumerate}
\item use the local Torelli theorem to prove that each K3 surface can be deformed to a K3 surface $X_0$ with the $Pic(X_0)$ generated by $L $ with $L^2=4$.
\item then to show that $X_0$ is isomorphic with a quartic surface.
\end{enumerate}\label{LP}
For \HK{} fourfolds that are numerical $K3^{[2]}$ similar steps were followed by O'Grady in order to prove Conjecture \ref{conj}:

 First using the local Torelli theorem the following was proved in \cite[Prop.~3.2, Prop.~4.1]{O}:
\begin{thm}[O'Grady]\label{T} Let $M$ be a symplectic $4$-fold numerically equivalent to
$(K3)^{[2]}$. There exists an IHS manifold $X$
deformation equivalent to $M$ such that:\begin{enumerate}
           \item $X$ has an ample divisor $H$ with $(h, h) = 2$ (i.e.~ $H^2=12$), where $h:=c_1(H)$,

\item $H^{1,1}_\mathbb{Z} (X) = \mathbb{Z}h$,
\item if $ \Sigma \in Z_1(X)$ is an integral algebraic 1-cycle on X and $cl(\Sigma) \in H_{\mathbb{Q}}(X)$
    is its Poincar\'e dual, then $cl(\Sigma) = mh^3 /6$ for some $m \in \mathbb{Z}$,
 \item if $H_1 , H_2 \in |H|$ are distinct then $H_1 \cap H_2$ is a reduced irreducible surface,
\item if $H_1 , H_2 ,H_3 \in |H|$ are linearly independent, the subscheme $H_1 \cap H_2 \cap H_3$
    has pure dimension $1$ and the Poincar\'{e} dual of the fundamental cycle $[H_1 \cap H_2 \cap H_3]$ is equal to $h^3$,
\item $\chi(\mathcal{O}_X(nH))=\frac{1}{2}n^4+\frac{5}{2}n^2+3$, $n \in \mathbb{Z}$.
\end{enumerate}
          \end{thm}
Thus in order to complete the proof of Conjecture \ref{conj} the problem is the second step. Such that we have to understand the geometry of \HK{} fourfolds with a polarisation of degree $q=2$.          

 Let us fix $X$ and $h:=c_1(H)$ as in Theorem \ref{T} above. By Theorem \ref{T}(6) and the Kodaira vanishing theorem we
infer $h^0(\mathcal{O}_X(H))=6$.  O'Grady consider the map given by the complete linear system $|H|$:       
  $$\varphi_{|H|}:X\dashrightarrow X'\subset \mathbb{P}^5$$ and prove in \cite{O} that it is either birational
  with a hypersurface of degree $6\leq d\leq 12$ or $2:1$ to a sextic hypersurface. In the second case $X'\subset\PP^5$ have to be an EPW sextic and in consequence $X$ of type $K3^{[2]}$.  So in order to prove Conjecture \ref{conj} it is enough to show that $\varphi_{|H|}$ cannot be birational. Our main result is the following contribution to the O'Grady program:
\begin{thm}\label{main}

 If the linear system $|H|$ on $X$ (as in Theorem \ref{T}) defines a birational map $$\varphi_{|H|}:X\dashrightarrow X'\subset \mathbb{P}^5$$ onto
 its image then $|H|$ has $0$-dimensional base locus of length $l\leq 3$. Moreover, the degree of $X'\subset \PP^5$ is $9\leq d \leq 12$.
\end{thm}
%Note, that if the base locus $Bs(|H|)$ of $|H|$ is finite then $\deg X'=12-l$. In particular the Theorem implies that the degree of $X'\subset\PP^5$ is in the range  $9\leq d\leq 12$.
In Section \ref{last} we give evidence that our proof of Theorem \ref{main} works also in the case $d=9$: thus we expect that the degree of $X'\subset \PP^5$ is $10\leq d\leq 12$.

It follows from Theorem \ref{main} that in order to complete the proof of the Conjecture \ref{conj} it is enough to show that $\varphi_{|H|}$ cannot be birational with image being a hypersurface in $\PP^5$ of degree $9\leq d\leq 12$.
Our second aim of this work (see Section \ref{rem}) is to study the geometric properties of the image $X'\subset \PP^5$ in the case when $\varphi_{|H|}$ is birational and $9\leq d \leq 12$.
% (that could lead to a possible counterexample to Conjecture \ref{conj}).

The idea of proof of Theorem \ref{main} is to argue by contradiction.
We assume that $\varphi_{|H|}$ defines a birational map with a hypersurface $X'\subset \PP^5$ of degree $6\leq d\leq 12$. In order to prove Theorem \ref{main} we have to obtain a contradiction in the following cases: 

\begin{itemize}
\item when the degree $d=6$; the contradiction is obtained in Section \ref{degree6},
\item when the degree $7\leq d\leq 8$;
%and the image $X'\subset \PP^5$ is non-normal such that its singular locus admits a sub-scheme structure
%$\mathcal{C}\subset X'\subset \PP^5$ defined by the conductor of the normalization. 
the contradiction is obtained 
%from the fact that the threefold $\mathcal{C}$ with the required properties cannot exist see 
in Section \ref{d89},
\item and when $9\leq d\leq 12$ and the base locus of $|H|$ is $\geq 1$ dimensional, we obtain a contradiction in Section \ref{rem}.

\end{itemize}
In order to explain the method of the proof let us discuss more precisely the second case.

 First, if the image $X'\subset \PP^5$ have degree $7$ or $8$ it is non-normal and his singular locus admits a sub-scheme structure
$$\mathcal{C}\subset X'\subset \PP^5$$ defined by the conductor of the normalization. 
In order to describe the threefold $\mathcal{C}\subset \PP^5$
we study $h^i(\mathcal{I}_{\mathcal{C}|\PP^5}(k))$ for $i,k\in \mathbb{Z}$. The simplest way to do this is to apply the projection formula that works for finite morphisms to relate those cohomologies to the cohomologies $h^i(\oo_{X}(n))$.
The problem is that the map $$\varphi_{|H|} \colon X\dashrightarrow X'\subset \PP^5$$ can contract some curves.
That is why we choose a generic codimension $2$ linear section
$X_D'\subset \PP^3$ of $X'\subset \mathbb{P}^5$ that avoids the image of the
curves contracted by $\varphi_{|H|}$ (see Proposition \ref{C}). Denote by $D$ the pre-image of
$X_D'$ on $X$. We construct the normalization $Y_D$ of $X_D'$ by a
sequence of blow-ups an blow-downs of $D$. Since the surface $Y_D$ has
rational singularities and the normalization $$Y_D\to X'_D$$ is finite, we can now find the possible cohomology
tables of the ideal of curves $C$ being a codimension $2$ linear section of $\mathcal{C}$. The main problem here is the explicit construction of the normalization $Y_D$. Then the cohomologies are computed in a standard way explained in detail in Section \ref{123}.
Next we address the problem of existence of curves with given cohomology tables.
The method to study such curves $C\subset \PP^3$ is to use the liaison theory described in \cite{PS}, \cite{MP} and \cite{M}.
In Section \ref{liaison} we introduce tools from liaison theory needed in this paper.
We study $C\subset \PP^3$ case by case according to the degree $d$ of $ X'\subset \PP^5$. 
%hypersurface $X'\subset \PP^5$. 

The above methods do not work in the cases $d=11, 12$ because the curve $C\subset \PP^3$ with the given cohomology table exists.
In Section \ref{last} we discuss the case $d=11$ where moreover no curves are contracted by $\varphi_{|H|}$. We show that the possible counterexample to Conjecture \ref{conj} in this case $d=11$ are related to a very special quintic hyprsurface in $\PP^5$ that we plan to study in a future work. The case $d=12$ is strongly related with EPW sextics and is discussed in \cite{K} where we give evidence that this case cannot occurs.

\section*{Acknowledgements} The author would like to thank Kieran O'Grady and Laurent Gruson for many valuable comments.
He also thanks Christian V. Bothmer, S{\l}awomir Cynk, Micha{\l} Kapustka, Christian Okonek, Szymon
Pli\'{s}, and Frank-Olaf Schreyer for helpful discussions.
\section{Preliminaries}\label{preliminary}
\subsection{The singular locus of $X'\subset \PP^5$}\label{cond} Our method to analyze non-normal hypersurfaces is to study their singular locuses.
Let us recall how to introduce a scheme structure on the singular set.

For a hypersurface $X'\subset \mathbb{P}^5$ we consider the ideal defined by the conductor.
 Let $$\beta:Y\rightarrow
X'\subset \mathbb{P}^5$$ be the normalization.
Set
$\mathrm{c}:=\mathrm{Ann}(\beta_{\ast}\mathcal{O}_Y/\mathcal{O}_{X'})\subset
\mathcal{O}_{X'}.$ 
Since $X'\subset \PP^5$ is a hypersurface,
 $\omega_{X'}$ is invertible. From the fact that $X'$ satisfies the Serre condition $S_2$,
we have
$$\mathrm{c}=\mathrm{Hom}_{\mathcal{O}_{X'}}(\beta_{\ast}\mathcal{O}_Y,\mathcal{O}_{X'})$$
(see \cite[p.~703]{R}). Assuming $Y$ is Cohen--Macaulay, we have
\begin{equation}\label{2}\mathrm{c}=\beta_{\ast}(\omega_{Y})\otimes_{\mathcal{O}_{X'}}\omega_{X'}^{-1}\end{equation}
(see \cite[p.~26]{S}).
The following important result is proved in \cite[p.~60]{Z}:
\begin{thm}[Zariski]\label{Za} Given that $X'\subset \PP^5$ is a hypersurface of dimension $r$, the conductor
ideal $\mathrm{c}$ is an unmixed ideal of dimension $r-1$ in
$\mathcal{O}_{X'}.$\end{thm} Denote by $\mathcal{C}\subset X'
\subset\mathbb{P}^5$ the subscheme of pure dimension $r-1$ defined
by $\mathrm{c}$. The codimension $2$ linear section of $\mathcal{C}$ is a locally Cohen--Macaulay curve $C\subset \mathbb{P}^3$.
We shall study the liaison class of this curve.

\subsection{Liaison theory}\label{liaison}
%The method of studying the threefold being the singular locus $\mathcal{C}\subset X'\subset \PP^5$ is to consider his curve section and study the properties of this curve.
%We need for this the Liaison theory. 
Let us thus recall some basic results from the
liaison theory \cite{MP,PS} (see \cite{GLM}) that we use in order to study $C\subset\PP^3$.

Let $C,D\subset\mathbb{P}^3$ be two locally Cohen--Macaulay curves that are
(algebraically) linked (see \cite[Def.~1.1]{MP}) through a complete intersection $X$ of
surface $s$ and $t$; we say then \emph{$s\times t$ linked}. Then $\deg C+\deg
D=st$ and
\begin{equation}\label{leq1}h^0(\mathcal{I}_C(n))-
h^0(\mathcal{I}_X(n))=h^2(\mathcal{I}_D(s+t-4-n))\end{equation}
\begin{equation}\label{leq2} h^1(\mathcal{I}_C(n))=h^1(\mathcal{I}_D(s+t-4-n))\end{equation}

For a curve $C\subset\mathbb{P}^3$ denote by
$M_C:=\oplus_{n\in\mathbb{Z}}H^1(\mathcal{I}_C(n))$ its Hartshorne--Rao module.
It is known that this gives a bijective correspondence between the
even liaison classes of such curves and graded $S:=k[x,y,z,t]$
modules of finite length, after identifying those that differ only
by shift. The curves with the smallest degree in a liaison
class are called minimal.

In \cite{MP} the authors show how to construct the minimal curves
from $$0\to L_4\to L_3\to L_2 \overset{\rho_2}\to L_1\to L_0\to
M_C\to 0$$ the minimal free resolution of its Hartshorne--Rao module. Note that
the minimal resolution of the curve $C$ is then as follows
$$0\to L_4\to L_3\oplus \bigoplus S(-r_j)\to \bigoplus S(-s_i)\to I_C\to 0 ,$$
where $r_i,s_i\in\mathbb{N}$.

More precisely they showed that, given a graded module $M$ of finite
length with $L_2=\bigoplus_{i=1}^{k}S(-b_i)$, there exists a
function $q \colon \{b_l,..., b_k\}\to \mathbb{N}$ (constructed using the
ranks of some submatrix of $\rho_2$ see \cite[p.287]{GLM}) such that
$B$ is a minimal curve for the biliaison class associated to the
module $M$ if and only if the homogeneous ideal $I_B$ is as follows:
$$0\to \mathcal{E} \to \bigoplus S(-n)^{l_2(n)-q(n)}\to I_B(h)\to 0.$$
Then the minimal curves $B'$ from the odd liaison class has
resolution
$$0\to \mathcal{N} \to S(a_0)\oplus S(a_1)\oplus \bigoplus S(n)^{q(n)}\to I_{B'}(h+a_0+a_1)\to 0,$$
where $\mathcal{E}, \mathcal{N}$ are modules described in \cite[IV]{MP} and
$$ h+\deg L_1- \deg L_0=\sum_{n \in \mathbb{Z}} n \cdot q(n)$$
is the shift of the Hartshorne--Rao module of the minimal curve $B$. Moreover
$a_0$, $a_1$ are constants introduced in \cite[Def.~2.4IV]{MP} such
that $B$ is $a_1+h\times a_0+h$ linked to a minimal curve from the
odd liaison class. We can also deduce the degrees and genus of
these curves as is shown in \cite[p.287]{GLM}. 
%Our aim is thus to find the $q$ functions of the liaison classes of curves defined by
%the conductor ideal in $\mathbb{P}^3$. We obtain in particular that
%the required Hartshorne--Rao module does not exists if $d=7,8$.

We shall also study other useful invariants. Recall from \cite{S2}, \cite{Sl} that $h_C(n)=\delta^2
h^2(\mathcal{I}_C(n))$ is called the spectrum of $C$. Then
$h_C(n)\geq 0$, $$\deg C=\sum_{n\in\mathbb{Z}}h_C(n),$$ and
$$p_a(C)=\sum_{n\in\mathbb{Z}}(n-1)h_C(n)+1.$$
where given a numerical function $f$, we let $\delta f$ denote
its first difference function $f(n)- f(n - 1)$.
If $C$ is obtained from $D$ by an elementary biliaison of height
$h=1$ on a surface of degree $s$ (see \cite[Def.~2.1 III]{MP}) then
\begin{equation}\label{leq3} h^0(\mathcal{I}_C(n))=
h^0(\mathcal{I}_D(n-1))+{n-s+2\choose 2}\end{equation}
\begin{equation}\label{leqh} h_C(n)=h_D(n-1)+h_F(n)\end{equation}
where $F$ is a plane curve of degree $s$.

\subsection{The reduction to codimension $2$ linear section of $X'\subset \PP^5$}\label{technical}

The aim of this section is to introduce technical results concerning 
the restriction of $\varphi_{|H|}$, for $H$ as in Theorem \ref{T}, to a codimension $2$ subvariety in order to obtain a finite morphism.

Consider the diagram
\begin{diagram}\label{d1}
X &\rDotsto^{\varphi_{|H|}} &X'\\
\uTo^{\pi}  & &\uTo^{\beta}\\
 {\overline{X}} & \rTo^{\rho}       &{Y},
\end{diagram}
where $\pi$ is a composition of blow-ups
with smooth centers such that the strict transform of the linear system $|H|$ is base point free. We will call the manifold
$\overline{X}$ a \textit{Hironaka model} of $(X,H)$ 
(after \cite[(1.4)]{F}, \cite{H}). The composition $\beta \circ \rho $ is the
Stein factorization of the birational morphism
$\overline{X}\rightarrow X$ induced by $\varphi_{|H|}$.
%The following proposition is one of the central point of this paper:
\begin{prop}\label{C} Let $H$ be as in Theorem \ref{T} and assume that the image $\beta(\rho(E))$, where $E$ is the exceptional locus of $\pi$,
is the sum of a finite number of linear subspaces of $\mathbb{P}^5$ with dimensions $\leq 3$. Then there exist two independent divisors $H_1,H_2\in |H|$ such that the restriction $\varphi_{|H|}|_{H_1\cap H_2}$ does not contract curves on $X$. 
\end{prop}

The proof the Proposition \ref{C} immediately follows from the following two Lemmas.% \ref{lemm3} and \ref{lemm4}:
\begin{lem}\label{lemm3} The map $\varphi_{|H|}$ does not contract surfaces on $X$ to points.
\end{lem}
\begin{proof} Suppose that $\varphi_{|H|}$ contracts a surface $S$ on $X$ to a
point $P\in\mathbb{P}^5$. Let us choose two independent hyperplanes
in $\mathbb{P}^5$ passing through $P$. It follows from
Theorem \ref{T} that their intersection is an irreducible
surface, which is a contradiction since $S$ is its proper component.
\end{proof}
\begin{lem}\label{lemm4} Assume that the image $\beta(\rho(E))$, where $E$ is the exceptional locus of $\pi$,
is the sum of a finite number of linear subspaces of $\mathbb{P}^5$ with dimensions $\leq 3$.
Then the map $\varphi_{|H|}$ does not contract divisors on $X$ to smaller dimensional subschemes.
\end{lem}
\begin{proof} Suppose that an irreducible divisor $D$ is contracted to a surface $S\subset \mathbb{P}^5$.
From Theorem \ref{T}(2) it follows that there exists $k\in
\mathbb{Z}$ such that $D\in |kH|$ with $k\geq 1$.

We claim that the surface $S$ is contained in $\beta(\rho(E))$.
Indeed suppose there is a curve $C$ that contracts to a point
outside $\beta(\rho(E))$. Then $C$ is disjoint with the base locus.
Since $H$ is ample we have $C\cdot H>0$ thus the pre-image of $C$ on
$\overline{X}$ cannot be contracted. The contradiction proves the
claim.

It follows that $S$ in contained in the sum of linear subspaces
$F_1,\dots,F_s$ of $\mathbb{P}^5$ with dimensions $\leq 3$.
Let us consider two cases.

If $s=1$ choose a generic hyperplane $\mathbb{P}^5\supset R_1 $ that contains $F_1$.
Then since $S\subset R_1$ the divisor $H_1\in |H|$ corresponding to $R_1$ contains $D$ as a proper component,
this is a contradiction with $D\in |kH|$ with $k\geq 1$.

Suppose that $s>1$, then since $S$ is irreducible (because $D$ is
irreducible), we deduce that $S$ is contained in one of the linear
space $F_1,\dots, F_s$. We obtain a contradiction as
before.\end{proof}
\begin{rem}\label{r3}
 O'Grady observed that the statement of Lemma \ref{lemm4} is equivalent to the following: the set
$$\{ p\in X' | \dim\varphi_{|H|}^{-1}(p)\geq 1 \}$$ has dimension at most 1 (here $\varphi_{|H|}^{-1}(p)$ is the set of  $x\in X$ outside the base-scheme such that $\varphi_{|H|}(x)=p$).
\end{rem}
\subsection{The intersection of divisors from $|H|$}
Let $H$ be as in Theorem \ref{T}. The aim of this section is to prove technical results needed in the proof of our main Theorem \ref{main}.

Three generic independent elements $H_2,H_3,H_4\in |H|$ intersect
along a subscheme $S$ of pure dimension $1$.
 Denote by $[S]\in Z_d(X)$ the fundamental cycle
associated to $S$ (as in \cite[p.~15]{Fu}). There is a unique
decomposition
\begin{equation}\label{gh} [S]=\Gamma+\Sigma,
\end{equation}
 where $\Gamma$ and $\Sigma$ are effective
$1$-cycles such that $$\supp\Sigma\subset \supp B$$ where $B$ is the base locus of
$|H|$ and $\supp \Gamma$
intersect $\supp \Sigma$ in points. We have
$$12= \deg (H\cdot(\Gamma + \Sigma))$$ (see \cite[\S 2]{O}). From Theorem
\ref{T}(3), we infer that the Poincar\'{e} dual, $cl(\Sigma)$ equals
$m h^3/6$. 

\begin{lem}\label{4}If $H_1\in |H|$ is generic,
$$12= d+\sum_{p\in \supp
B}mult_p(H_1\cdot\Gamma)+2m.$$
\end{lem}
\begin{proof} In fact the above equation holds in the following
situation: let $\Theta\subset |H|$ be a $4$-dimensional linear
subsystem. Then Lemma \ref{4} holds for arbitrary linearly independent
$H_1,\dots, H_4\in\Theta $ if Eqnt.~(4.0.24) of \cite{O} holds for
one set of linearly independent $H_1,\dots, H_4\in \Theta$.
\end{proof}
Let us prove the following:
\begin{lem}\label{lemm2} If the base locus of $|H|$ is $0$-dimensional, then the
generic divisor in $|H|$ is smooth.
\end{lem}
\begin{proof}(cf.~\cite[(2.5)]{F})
Let $\pi_1:T\rightarrow X$ be a blow-up of a point from the base
locus such that $E$ is the exceptional divisor. Then
$\pi^{\ast}(H)-sE$ is semi-positive, thus
$$0\leq(\pi^{\ast}(H)-sE)^4=12-s^4,$$ so $s=1$.
\end{proof}
More generally, O'Grady observed the following:

\begin{lem}\label{lemmO}
If the base locus $B$ of $|H|$ is $0$-dimensional, then the
intersection $D$ of two generic elements of $|H|$ is smooth.
\end{lem}
\begin{proof}
By the Bertini theorem $D$ is smooth outside $B$. Let $\Theta\subset
|H|$ be a generic $4$-dimensional linear subspace and
$$m:=4-\dim \bigcap_{H'\in\Theta}T_{p_0}(H').$$
We may choose linearly independent $H_1,\dots,H_4\in\Theta$ such that
$$T_{p_0}(H_1)\cap \dots\cap T_{p_0}(H_4)=\bigcap_{H'\in\Theta}T_{p_0}(H').$$
Since $\dim B=0$ the intersection $H_1\cap\dots \cap H_4$ is proper and hence
$$12=d+\sum_{p\in \supp B}mult_P(H_1\dots H_4).$$
Now $d\geq 7$ and hence $\mult_P(H_1\dots H_4)\leq 5$ for all $p\in \supp B$, in particular for $p=p_0$.
It follows that $m\geq 2$ and hence the intersection of two generic divisors in $\Theta$ is smooth at $p_0$. Notice that if $d\geq 9$ we actually get $m\geq 3$.
\end{proof}
\begin{rem}\label{k}
 The above argument shows also the following: if $B=B^1\coprod Z$ where $\dim Z=0$ and $B^1$ is the sum of $1$-dimensional components, then there exists a surface containing $B$
 which is smooth at each point of $Z$ and hence the scheme $B$ is planar at each of its isolated points.
Moreover, it is curvilinear (contained in a smooth curve) if
$d\geq9$. It follows that in these cases that the base-point locus
$B$ of $|H|$ is of length $12-d$.\end{rem}

\begin{rem}\label{lp}
Let $\pi\colon \tilde{X}\rightarrow X$ the blow-up of the
base-scheme $B$; thus $\varphi_{|H|}$ defines a regular map
$\tilde{\varphi}\colon \tilde{X}\rightarrow \mathbb{P}^5$. Let $b\in
B$ be an isolated point. If $B$ is a local complete intersection at
$b$ then $\pi^{-1}(b)$ is irreducible of dimension $3$ and moreover
$\tilde{\varphi}(\pi^{-1}(b))$ is a $3$-dimensional linear subspace
of $\mathbb{P}^5$. In particular if $\dim B=0$ and $B$ is a local
complete intersection we get that the hypothesis of Lemma
\ref{lemm4} is satisfied and hence for generic $D$ there are no
contracted curves on $D$. Suppose that $d\geq 9$ and $\dim B=0$; by
Remark \ref{k} we get that $B$ is curvilinear, in particular a
l.c.i..\end{rem}
\subsection{The geometry of the base locus of $|H|$}
Let us again consider $H$ as in Theorem \ref{T}.
In this section we study general properties of the base locus of $|H|$.

Let $D\subset X$ be the intersection of two generic divisors from
$|H|$ and $X'_D\subset\mathbb{P}^3$ the corresponding linear section
of $X'\subset \PP^5$. From Theorem \ref{T}(4) the surface $D$ is reduced and
irreducible. Since from Proposition \ref{lemmN}
the one dimensional part of the base locus $B^1$ is reduced at the
generic point we infer by the Bertini theorem (cf.~\cite[Thm.~2.1]{DH}) that the
surface $D$ has isolated singularities. Thus by the Serre
criterion the surface $D$ it is normal. Moreover, $D$ is locally a complete intersection, so it
is locally Cohen--Macaulay and $\omega_D=2H|_D$. 

If $H_1,H_2,H_3\in
|H|$ are generic we write as above $$[H_1\cap H_2\cap
H_3]=\Gamma'+\Sigma'$$ where $\Sigma'$ and $\Gamma'$ are cycles such that $\supp\Sigma'\subset \supp B$ and $$\dim(\supp
\Gamma'\cap\supp \Sigma')=0.$$ The following diagram is induced by restriction from
diagram (\ref{d1}):
\begin{diagram}\label{dN}
D &\rDotsto^{|H'_D|} &X'_D\subset \mathbb{P}^3 \\
\uTo^{\pi_D}  & &\uTo^{\beta_D}\\
 {\overline{D}} & \rTo^{\rho_D}       &Y_D
\end{diagram}
here $H_D'$ is the restriction of $H$ to $D$. Denote by $H_D$ the
pull back of $\mathcal{O}_{X_D'}(1)$ by $\beta_D\circ \rho_D$.
One infer the following:
\begin{prop}\label{lemmN}
Suppose that $d\geq 7$. The base-scheme $B$ is reduced at the
generic point of any of its $1$-dimensional irreducible components.
\end{prop}
\begin{proof} We argue by contradiction. Suppose that the cycle $\Sigma$
associated to generic $$H_2,H_3,H_4\in |H|$$ is non-reduced. By
 Lemma \ref{4} one gets that $d=7$ and $\Sigma=2\Sigma'$ where
$\Sigma'$ is an irreducible curve. Moreover, $B$ is of pure dimension $1$,
non-reduced irreducible and there is a unique point $p\in
\Sigma\cap\Gamma$ such that \begin{equation}\label{w1} \mult_p(H_1\cdot
\Gamma)=1 \ \ \ \ H_1\in |H|\ \ \ \textrm{generic}. \end{equation}
We claim that for generic $H_2,H_3,H_4\in |H|$ we have
\begin{equation}\label{op}  (T_pH_2)\cap(T_pH_3)\cap(T_pH_4)= T_pB \ \ \ \ \ \ \forall p\in B. \end{equation}
First notice that \begin{equation}\label{oi} \dim T_p B\geq 2 \ \ \
\ \ \ \forall p\in B\end{equation}  because $B$ is everywhere
non-reduced and of dimension $1$. Moreover, $\dim T_p B=2$ for
generic $p\in B$ because $\Sigma=2\Sigma'$. Let $H_3,H_4\in |H|$ are
generic then
\begin{enumerate}
\item[(a)] $\dim((T_pH_3)\cap(T_pH_4))=2$, if $\dim((T_pH_4)=3$ (thus
for generic $p$),
\item[(b)] $\dim((T_pH_3)\cap(T_pH_4))=4$, only if
$\dim((T_pB)=4$.\end{enumerate} Now choose a generic $H_2\in|H|$;
then (\ref{op}) follows from (a) and (b) above and (\ref{oi}). Let
$$[H_2\cap H_3\cap H_4]=\Sigma+\Gamma,$$ as usual, and $p$ be the
unique point in $\Sigma \cap  \Gamma $. Let $H_4\in |H|$ be
arbitrary; then $T_p H_1\supset T_p B$ by (\ref{op}). It follows
that $$\mult_p(H_1\cdot \Gamma)\geq 2.$$ We obtained a contradiction
with (\ref{w1}).
\end{proof}
\begin{cor}\label{cor}Suppose that $d\geq7 $. Let $H_2,H_3,H_4\in |H|$ are
generic and $$[H_2\cap H_3\cap H_4]=\Sigma+\Gamma$$ where $\Sigma$ is
supported on $ B$ and $$\supp(\Gamma)\cap \supp(\Sigma)$$ is
$0$-dimensional. Then the cycle $\Sigma$ is reduced.
\end{cor}
\begin{proof} The surface $D=H_2\cap H_3$ is reduced, irreducible and normal. Let
$\overline{\Gamma}+\overline{\Sigma} $ be the Cartier divisor on $D$
corresponding to $H_4$ such that $\supp
\Sigma=\supp\overline{\Sigma}$ and $\supp
\Gamma=\supp\overline{\Gamma}$. It is enough to prove that the
scheme $\overline{\Sigma}$ is reduced at his generic point $q$.
Suppose the contrary, then $$T_q\overline{\Sigma}=T_q D$$ and this
holds for a generic $H_4$. This is a contradiction since $B$ is
reduced at $q$.
\end{proof}
\begin{rem}\label{poi} 
Let $D=H_1\cap H_2$ where $H_1,H_2\in
|H|$ are generic, such that $D$ is normal.
Suppose that for $$H_2,H_3,H_4\in |H|$$
generic $\Gamma'$ and $\Sigma'$ are Cartier divisors on $D$. Then the linear system $|\Gamma'|+\Sigma'$ can be
naturally identified with the linear system $|H||_D$ being the
restriction of $|H|$ to $D$. Indeed, it is enough
to observe that they have the same dimension. This follows from the
fact that $H^1(\mathcal{O}_{H_1}(H))=0$, since $H$ is ample. 
Thus
$|H||_D$ is a complete linear system. It follows from Corollary
\ref{cor} that the $1$-dimensional part of the base locus of
$|H||_D$ is reduced and equal to $\Sigma$.
%We conclude that $\Lambda_D$ has fixed part equal to the
%$1$-dimensional components of $B$ (the base scheme of $|H|$).
\end{rem}

\section{Proof of Theorem \ref{main} in the case $d=6$}\label{degree6} Let us first
consider the case $d=6$ of the Theorem \ref{main}. Let us prove more generally, that there does not exist a \HK{} manifold with a divisor giving a birational morphism to a normal hypersurface
of degree $6$. Suppose the contrary. 
Denote by $Z$ the subscheme of $\mathbb{P}^5$ defined by the adjoint ideal $\mathrm{adj}(X') \subset \mathcal{O}_{\mathbb{P}^5}$
(see \cite[Def.~9.3.47]{L}). From \cite[Prop.~9.3.48]{L} one has the
following exact sequence: \begin{equation}\label{1} 0\rightarrow
\mathcal{O}_{\mathbb{P}^5}(-6)\rightarrow
\mathcal{O}_{\mathbb{P}^5}(d-6)\otimes\mathcal{I}_Z\rightarrow
\beta_{\ast}(\mathcal{O}_{Y}(K_Y))\rightarrow 0.\end{equation} 
It follows that
$h^0(\mathcal{I}_Z)=1$, thus
$\mathcal{I}_Z=\mathcal{O}_{\mathbb{P}^5}$. From
\cite[Prop.~9.3.43]{L}, we infer that $X'\subset\mathbb{P}^5$ is a
normal hypersurface of degree $6$ that has
rational singularities. Let $(\overline{X},\overline{H})$ be a
Hironaka model of $(X,H)$. Then $|\overline{H}|$ gives a morphism
$\rho:\overline{X}\rightarrow X'$ that is a resolution of $X'$. So
from \cite[Thm.~5.10]{KM}, we infer
$R^i\rho_{\ast}(\mathcal{O}_{\overline{X}})=0$ for $i>0$. In
particular,
$$h^2(\mathcal{O}_{\overline{X}})=h^2(\rho_{\ast}(\mathcal{O}_{\overline{X}}))=h^2(\mathcal{O}_{X'})=0.$$
However, from the Hodge symmetry we infer
$$h^2(\mathcal{O}_{X})=h^0(\Omega^2_{X})=1.$$ Since the resolution
$\overline{X}\rightarrow X$ is obtained by a sequence of blow-ups
we obtain $h^2(\mathcal{O}_X)=h^2(\mathcal{O}_{\overline{X}}),$ a
contradiction.
\section{Proof of Theorem \ref{main} in the cases $d=7,8$}\label{d89}
Let $H$ be an ample divisor on $X$ with $q(H)=2$ as in Theorem \ref{T}.
%In order to complete the proof of theorem \ref{main} one have to prove the following statements:
%\begin{itemize}
% if the map $\varphi_{|H|}$ is birational then the linear system $|H|$ has a $0$ dimensional base locus,
The aim of this section is to show the following implication: if the map $\varphi_{|H|}$ is birational then the degree of the image $\varphi_{|H|}(X)=X'\subset \PP^5$ cannot be $7$ or $8$.
%\end{itemize}
%We shall treat those two cases separately in the next sections
%\subsection{Zero dimensional base locus}
%\section{Degrees 8 and 7}
%The aim of this section is to complete the proof of Theorem \ref{T1} and show that the degree $d$ of $X'\subset \PP^5$ cannot be equal to $7$ or $8$. 

In the range $d=7,8$ two situations are possible: either $|H|$ has $0$ dimensional base locus or a $1$-dimensional base locus.
We shall treat each of those cases in the next subsections arguing by contradiction.
\subsection{Suppose first that the base locus $B$ is $1$-dimensional.} If
$H_1,H_2,H_3\in |H|$ are generic we write as in Equation (\ref{gh}) $$H_1\cap H_2\cap
H_3=\Gamma'+\Sigma'.$$ Then the cycle $\Sigma'$ is not zero and
reduced. Denote as above by $D$ the intersection of $H_1$ and
$H_2$. Moreover, $B^1$ is the union of
$1$-dimensional components of the base scheme $B$. Let us first prove some more technical results that will distinguish three possible subcases.
\begin{lem}\label{clai}
 Let $\Theta \subset |H|$ be a generic $4$-dimensional linear subsystem i.e.~for which Lemma \ref{4} works.
 %\label{gh}
Given a generic $p\in \supp B^1$ there exist
linearly independent $H_2',H_3',H_4'\in \Theta$ such that $p\in
\supp \Gamma$, where $\Sigma+\Gamma =[H_2'\cap H_3'\cap H_4']$ 
%are as in (\ref{gh}). 
\end{lem}
\begin{proof}
% Let us consider the first statement.
 We can assume
that $p\in B^1$ is smooth. Let us show that there are three
independent elements of $\Theta$ such that their intersection is
singular at $p$. Since $\Theta$ is $4$-dimensional, we can find two
elements $H_2',H_3'$  in $\Theta$ that have the same tangent space at $p$.
Their intersection $D_1$ is singular at $p$ but from Theorem
\ref{T}(4) irreducible and reduced. It is enough to choose the
third element $H_4'$ generically i.e.~such that $H_4'$ cuts $D_1$
transversally at a generic point of $B^1$.
\end{proof}
\begin{lem}\label{clai1}  Let $\Theta \subset |H|$ be a generic $4$-dimensional linear subsystem.
Suppose that there exists $p_0\in \supp B^1$ such that one
has $p_0\in \supp \Gamma$ for
a generic set of linearly independent $H_2',H_3',H_4'\in \Theta$, where $$\Sigma+\Gamma =[H_2' \cap H_3' \cap
H_4']$$ is as in Equation (\ref{gh}). Then there is a unique such $p_0$, $d=7$ and for
a generic set of linearly independent $H_2',H_3',H_4'\in\Theta$ the
corresponding (by Equation (\ref{gh})) cycle $\Sigma$ is reduced and irreducible.
\end{lem}
\begin{proof}
%Let us prove the second statement. 
Let $\Theta=\mathbb{P}(W)$ where
$W\subset H^0(X,\mathcal{O}_X(H))$ is a $4$-dimensional sub
vector-space. Given $p\in B^1$ consider the differential map
\begin{center} $\delta_p\colon W\rightarrow \Omega_p X$,\ \ \ \ \ \
      $   \delta_p(\sigma)=d\sigma(p).$\end{center}
Let $K_p:=\ker (\delta_p).$ Let $p$ be a generic point of a
component of $B^1$; then $\dim K_p=1$ by Proposition \ref{lemmN} and
if $U\subset W$ is a generic $3$-dimensional subspace containing
$K_p$ then letting $$\Sigma+\Gamma=[H_2'\cap H_3'\cap H_4']$$ where
$\langle H_2',H_3',H_4'\rangle=\mathbb{P}(U)$,  from Lemma \ref{clai} we have $p\in \Gamma$.
%($\Sigma$ is reduced ). 

Now let $p_0$ be as is the statement of the Claim; then
$\dim K_{p_0}\geq 2$. It follows that the subset of
$Gr(3,W)=\mathbb{P}(W^{\vee})$ defined by
$$\{ U|\ \delta_{p_0}(U)\neq im(\delta_{p_0}) \}$$
is a linear subspace of dimension at most $1$. Since the set of
$U\in Gr(3,W)$ containing $K_p$ is a $2$-dimensional linear subspace
there exists $U_0\in Gr(3,W)$ containing $K_{p}$ which does not
belong to the set of the above equation and such that the
corresponding $\Sigma$ is reduced. Let $$\langle
H_2',H_3',H_4'\rangle=\mathbb{P}(U_0).$$ If $\langle
H_1',H_2',H_3',H_4'\rangle=\Theta$ then
$$\mult_p(H_1'\cdot\Gamma)\geq 1 \ \ \ \text{and} \ \ \ \ \ 
\mult_{p_0}(H_1'\cdot\Gamma)\geq 2$$ hence the Claim follows from
Lemma \ref{4}.
\end{proof}
Thus one of the following cases happen:
\begin{enumerate}
\item either the divisor $\Gamma'$ is Cartier on $D$ and define a base-point-free linear system, or
\item  the divisor $\Gamma'$ is Cartier $D$ and $|\Gamma'|$ has only isolated base points that are outside $\Sigma'$, or
\item we have $d=7$, there is a unique point $P_0\in \Sigma'$ such that $P_0\in \Gamma'$
for each $\Gamma'$ (for a generic choice of $H_3'$).
\end{enumerate}
We shall treat each case separately arguing by contradiction. In the next subsections we obtain a contradiction in each of the above three cases.
\subsubsection{Assume we are in case (1). }\label{op1} 
Our aim is first to construct the normalization of the surface $X_D'\subset \PP^3$ in order to compute the possible cohomology tables of the ideal sheaf of his singular curve. We obtain a contradiction by showing that there are no curves with the given cohomology table.

If we are in case (1), by \cite[Thm.~2.1]{DH}, the
surface $D$ has only isolated singularities at singular points of
$\Sigma'$. Let $(\tilde{D},\Sigma^{\circ})$ be a minimal resolution
of $(D,\Sigma)$, $\tilde{\Gamma}$ the pull-back of $\Gamma'$ on
$\tilde{D}$, and $\Sigma^{\circ}$ the strict transform of $\Sigma'$.

 \begin{prop}\label{S} The morphism from $D$ given by
$|\Gamma'|$ is the normalization of $X'_D$.\end{prop}
\begin{proof} Let us first consider some technical results.
\begin{claim} $\Sigma'$ is irreducible. 
\end{claim}
%\begin{proof}
Indeed, if $\Sigma'$ has two
components $\Sigma_1$ and $\Sigma_2$. Denote by $\Sigma_1^{\circ}$,
$\Sigma_2^{\circ}$ the corresponding components of $\Sigma^{\circ}$.
Then from Lemmas \ref{clai} and \ref{clai1} we have $$\Sigma_1^{\circ}\cdot
\tilde{\Gamma}\geq 1\ \ \ \text{and} \ \ \ \Sigma_2^{\circ}\cdot \tilde{\Gamma}\geq
1$$ on $D$. This contradicts Lemma \ref{4}. The claim follows.
%\end{proof}

Observe that if $p\in \Sigma'$ is smooth then
$$\mult_p(\Sigma'\cdot \Gamma')=\mult_p(H_3\cdot \Gamma'),$$ where
$H_3\in |H|$ is generic.
Thus by Lemma \ref{4} we have 
\begin{equation}\label{gar}
\Sigma'\cdot \Gamma'= 2 \ \ \  \text{(resp.~3) if} \ \ \  d=8 \ \ \ 
\text{(resp.~7)}. 
\end{equation}
So, the image of $\Sigma^{\circ}$ by
$|\tilde{\Gamma}|$ is a smooth conic or a line in $X'_D\subset
\mathbb{P}^3$ (resp.~a rational normal curve in $\mathbb{P}^3$, a
smooth elliptic curve, a line, or a singular cubic curve in
$\mathbb{P}^2$).

 We have to prove that $\varphi_{|H|}$ does not contract
divisors on $X$ to surfaces.
% in order to cheek that the assumptions of Lemma \ref{lemm4} are filled.
Suppose that an irreducible divisor
$S\subset X$ is contracted to a surface. Let
$$\pi:\tilde{X}\rightarrow X$$ be the blow up of $\Sigma'\subset X$.
Then the image
$$I=\varphi_{|\pi^{\ast}(H)-E|}(E)\subset \mathbb{P}^5$$ is covered by
planes since $E$ is a $\PP^2$ fibration over $\Sigma'$.

First, assume that $\Sigma'\subset D$ maps by $|\Gamma'|$ to a line.
It follows that the generic $3$-dimensional linear section of $I$ is
a contained in a line. Thus, $I$ is contained in a $3$-dimensional
linear space in $\mathbb{P}^5$. We conclude by Lemma \ref{lemm4}.

 If the image of $\Sigma^{\circ}$ is not a line then $|\tilde{\Gamma}|$
is generically $1:1$ on $\Sigma^{\circ}$. We have two possibilities:
Either there is a point $P\in \Sigma'$
such that for a generic curve $C\subset S$ contracted by $\varphi_{|H|}$
we have $P\in C$ or there is no such point. In the first case the image of $C$ is contained in the plane in
$\mathbb{P}^5$ being the image of the fiber of the exceptional
divisor $E$ under $P$. We can conclude by Lemma \ref{lemm4}.

Or let us consider the remaining case. Suppose that $C\subset S$ is a
curve on $D$ contracted by $\varphi_{|H|}$. Denote by $C^{\circ}$ its
strict transform on $\tilde{D}$. We can assume that $C$ intersects
$\Sigma'$ in smooth points on $\Sigma'$ thus smooth on $D$. The
divisors $\tilde{\Gamma}|_{\Sigma^{\circ}}$ gives a linear system
$\Lambda$ on $\Sigma^{\circ}$. We have $$C^{\circ}\cdot
\Sigma^{\circ}=C\cdot H\geq 2.$$ Thus if $P\in\Sigma^{\circ}\cap
C^{\circ}$ and $P+A\in \Lambda$, where $A$ is an effective divisor
on $\Sigma^{\circ}$, then $$\supp A\cap C^{\circ}\cap \Sigma^{\circ}$$
is non-empty. Now, observe that the only linear systems of degree
$2$ with this property are one-dimensional. Moreover, the images of
such linear systems of degree $3$ that are not lines are singular,
so they are plane cubics.
% (the degree $\leq 3$ by Equation (\ref{gar})).

Let us assume that the image by $|\Gamma'|$ of $\Sigma'\subset D$ is
a plane cubic. It follows that the image of $E$ in $\mathbb{P}^5$ is
contained in a hyperplane $L$. Since $S\in |kH|$ with $k\geq 1$ and
the image of $S$ is contained in $L$ and cannot be a proper component of
the pre-image of $L$, we obtain $k=1$ and
$$C+\Sigma'\in |H||_D.$$ Thus $$C^{\circ}\cdot \Sigma^{\circ}=3,$$ so
the linear system $\Lambda$ has the property; if
$P\in\Sigma^{\circ}\cap C^{\circ}$ and $P+A\in \Lambda$, where $A$
is an effective divisor on $\Sigma^{\circ}$, then $$P+ A=
C^{\circ}\cdot \Sigma^{\circ}.$$ It follows that the image of
$\Sigma'$ is a line, a contradiction.
\end{proof}

We deduce \cite[Thm.~3.5]{KLU} that the conductor of the normalization of $X'_D$ defines
a locally CM subscheme $C\subset \mathbb{P}^3$ such that \begin{equation}\label{po}2\deg
C=\Gamma'((d-6)\Gamma'-2\Sigma').\end{equation} 
The above calculation will be explained in details in Section \ref{123}.
We obtain a contradiction if
$d=7$ because the left hand side of Equation (\ref{po}) is even and the right odd. 

So assume that $d=8$, thus $\deg(C)=6$.  In order to find the cohomology of the ideal sheaf of $C\subset \PP^3$ we use the following:
\begin{lem}\label{kp}
We have
$$h^0(\mathcal{I}_C(n+4))=h^0(K_D+n\Gamma')$$ for $n=-2,-1,0,1$. 
\end{lem}
We
infer $$h^0(K_D)=h^0(2\Gamma'+2\Sigma')=h^0(D,\mathcal{O}_D(2H)),$$
thus from Theorem \ref{T} (6) we obtain $h^0(\mathcal{I}_C(4))=11$.
Since $$h^0(\mathcal{I}_C(1))=0\ \ \ \ \ \text{and} \ \ \ h^0(2\Sigma')\geq 1,$$ we infer
$h^0(\mathcal{I}_C(2))=1$. It follows that
$h^0(\mathcal{I}_C(3))\geq 4$. Moreover, by the Riemann-Roch theorem
we find that $$\chi(K_D+n\Gamma')=4n^2+10n+12$$ and from the Castelnuovo Mumford criterium that the ideal of $C\subset \PP^3$ is generated by quintics. We deduce the following table:
\begin{center}
            %\begin{table}
            \renewcommand*{\arraystretch}{1}
\begin{table}[htp]
 ×

\begin{tabular}{|c|c|c|c|c|}

\hline
            $n$& $h^0(K_D+n\Gamma')$& $h^1(K_D+n\Gamma')$& $h^2(K_D+n\Gamma')$& $\chi(K_D+n\Gamma')$                       \\ \hline
                                                                   \hline

            $0$&11&$0 $&$1$& $12$           \\ \hline
            $-1$&$4\leq$ &$x$&$4$& $6$ \\ \hline
            $-2$& 1&$y$&$7+y$&$8$                    \\ \hline
            $-3$&0&$z$&$18+z$&$18$
            \\ \hline
            $-4$&0&$t$&$36+t$&$36$
            \\ \hline
            \end{tabular}
 \end{table}
           %\end{table}
            \end{center}
Here we have $4\geq y\geq 3$, $z\geq 2$, and $x\geq 2$. Let
$B\subset \mathbb{P}^3$ be a degree $4$ curve $2\times 5$ linked to
$C$. We infer from Lemma \ref{kp} the following:
\begin{center}
            %\begin{table}
            \renewcommand*{\arraystretch}{1}
\begin{table}[htp]
 ×

\begin{tabular}{c|ccc}

            $n$& $h^0(\mathcal{I}_B(n))$& $h^1(\mathcal{I}_B(n))$& $h_B(n)$                       \\ \hline

            $3$&$t+5$ &$t$&$0$ \\
            $2$&$z-1$ &$z$&$x-2$ \\
            $1$& $y-3$&$y$&$a$                    \\
            $0$&0&$x$&$b$

             \end{tabular}
 \end{table}
           %\end{table}
            \end{center}
            see Section \ref{liaison} for the definition of $h_B(n)$.
If $x>2$ then $h_C(2)\leq -1$ contradiction, thus $x=2$, $a=1$, and $b=3$.
It follows that $B$ is not extremal (see \cite{S3}) and that
$p_a(B)\geq-2$. We have the following inequalities from \cite{N} or
\cite[Thm.~4.4]{S3}: 
$$y\leq 3, \ \ \ \  z\leq 2, \ \ \  t\leq 1.$$ Thus we have
two possibilities $(y,z,t)=(3,2,0)$ or $(3,2,1)$. It follows that
$B$ is contained in a quadric and $$(h_{C}(0),h_{C}(1),h_{C}(2))$$ is equal to $(1,4,1)$ or $(2,2,2)$. 

We infer from \cite[Cor.~4.4]{S2} and Equation (\ref{leqh})
that $B$ is minimal in its biliaison class, and $C$ can be bilinked down on the quadric to $C_0$ a minimal curve of degree $2$. We obtain a
contradiction with \cite[Ex.~1.5.11]{M} where all the possible
deficiency modules of non reduced curves of degree $2$ are described.

\subsubsection{Assume we are in case (2).} 

Then the $0$-dimensional
components of the base locus of $|H|$ have length $\leq 2$
 and from Lemma \ref{4} the cycle $\Sigma'$ is reduced and irreducible. Thus from Lemma \ref{lemmO} the surface $D$ is
smooth outside $\Sigma'$, so from
\cite[Thm.~2.1]{DH} has only isolated singularities. 
\begin{lem} The composition $\rho_D\circ\beta_D$ gives the normalization of
$X'_D\subset \PP^3$.
\end{lem}
\begin{proof}
We have $\Gamma'\cdot \Sigma'\leq
2.$ So the image of $\Sigma'$ is a smooth conic or a line and we can
prove as in Proposition \ref{S} that $\varphi_{|H|}$ does not contract
curves on $D$. 
\end{proof}
So the conductor of the normalization of $X'_D$ defines a
locally CM subscheme $C\subset \mathbb{P}^3$ such that $$2\deg
C=\Gamma''((d-6)\Gamma''-2\Sigma'-R),$$ where $R$ is an effective
divisor supported on the exceptional lines on $\overline{D}$ and
$\Gamma''$ (resp.~$\Sigma''$) the strict transform of $\Gamma'$
(resp.~$\Sigma'$) on $\overline{D}$. Now, if $d=7$ we obtain a
contradiction with $\deg (C)\geq 1$.

Assume that $d=8$. Then $$\Gamma'\cdot \Sigma'=1$$ and $\Gamma'$ has
exactly one isolated simple base point $P_0$, denote by $E$ the
exceptional divisor of the blow-up at $P_0$. From the adjunction
formula we infer
$$2g(\Gamma'')-2=\Gamma''(\Gamma''+K_{\overline{D}})=\Gamma''(3\Gamma''+2\Sigma''+3E)=29,$$
a contradiction since the genus $g(\Gamma'')$ is an integer.

\subsubsection{Assume we are in case (3).}
\begin{lem}\label{lemmW} Under the assumption of case (3)
 the intersection $D$ of two generic
divisors $H_1',H_2'\in|H|$ is smooth at $P_0$. \end{lem}
\begin{proof} Arguing as in the proof of
\cite[Prop.~5.4(2)]{O} we infer that the generic $\Gamma'$ is smooth
at $P_0$; moreover, the tangent direction of $\Gamma'$ is not
contained in the tangent space $T_{P_0}\Sigma'$. If $H_1'$ is
singular at $P_0$, then the multiplicity of the intersection of
three generic divisors from $\Theta$ at $P_0$ is $\geq 8$. It
follows that $\Sigma'$ is singular and $T_{P_0}\Sigma'$ has
dimension $\geq 3$. Thus the tangent space $T_{P_0}\Sigma'$ cannot
intersect transversally $T_{P_0}\Gamma'$, a contradiction. Repeating
the above arguments for $H_1'$ instead of $X$ we end the proof.
\end{proof}
From
Lemma \ref{4} we infer that $P_0$ is a simple base point of $|\Gamma'|$,
i.e.~the strict transform $\Gamma''$ of $\Gamma'$ on the blowing-up
$D'$ of $D$ at $P_0$ is base-point-free.
\begin{lem} The morphism from $D'$ given by $|\Gamma''|$ does not contract
curves and gives the normalization
of $X'_D$.
\end{lem} 
\begin{proof}
It follows that $\Sigma'$ and $\Gamma'$ are Cartier divisors. Denote by $E $ the
exceptional divisor and by $\Sigma''$ the strict transform of
$\Sigma'$. It follows that we can resolve the indeterminacy of $|H|$
by blowing-up $\Sigma'$ and then the fiber over $P_0$ of the
obtained exceptional divisor. Denote by $\overline{X}$ the obtained
threefold and by $E_1\subset X$, $E_2\subset X$ the resulting
exceptional divisors.

\begin{claim} The morphism $\overline{\varphi}\colon
\overline{X}\rightarrow X'$ induced from $\varphi_{|H|}$ maps $E_1$ and
$E_2$ into two $3$-dimensional linear subspaces of $\mathbb{P}^5$.
\end{claim}
Indeed, it is enough to observe that $E$ and $\Sigma''$ maps into
generic hyperplane sections of $E_1$ and $E_2$. Now, the image of
$E$ is a line since $\Gamma'$ is smooth at $P_0$ and $P_0$ is a
simple base point. Since
$\Sigma''\cdot \Gamma''=1$ the image of $\Sigma'$ is also a line, the claim follows.

So, we can use Lemma \ref{lemm4} to conclude.
\end{proof}
 We infer that the conductor of this normalization
defines an CM subscheme $C\subset \mathbb{P}^3$ such that $$2\deg
(C)=\Gamma''(\Gamma''-2\Sigma''-sE)$$ such that $s\geq 5$. This is a
contradiction since $\Gamma''\cdot E=1$, $\Gamma''\cdot \Sigma''=1$,
and $\deg(C)>0$.

\subsection{Suppose next that the base locus $B$ is $0$-dimensional.}
As before, we denote by $D$ the intersection of two generic elements
of $|H|$ and set $H|_D=\Gamma'$. From Lemma \ref{lemmO} we infer
that $D$ is smooth. Let us consider cases are when $\supp B$ is one point
$P$, where $B$ is the base locus of $|H|$, the other cases are treated analogously. Let us consider the cases $d=7$ and $d=8$ separately.

\subsubsection{Suppose first that $d=8$.} As usual let us first construct a Hironaka model of $(D,\Gamma')$.
We assume that $\supp B=\{P\}$ of the base locus is one point. Note that the other cases are simpler and we obtain the same results. 
We consider two possibilities either $\Gamma'$ is smooth or singular at $P$.

If the generic $\Gamma'$ is smooth at $P$ 
we see that a Hironaka model
$\overline{D}$ is obtained by four blowings-up at each step of the
unique fixed point of the linear system $|\overline{\Gamma}|$ which
is the strict transform of the linear system $|\Gamma'|$. We have
however five possible configurations of the resulting exceptional
curve depending on the positions of choosen points on exceptional
divisors. The morphism $$\rho_D:\overline{D}\rightarrow
Y_D$$ is birational and contracts all the exceptional divisors except
the last one. From \cite[Thm.~3]{A} we infer that the singularities
of $Y_D$ are Du Val or rational triple points (see
\cite[p.~135]{A}). 
%We deduce the following:
%\begin{lem} The map $\rho_D$ does not
%contracts curves. 
%\end{lem}
%Since a surface with rational singularities is
%$\mathbb{Q}$-factorial and Cohen--Macaulay.
%we can argue as in the case $d=9$ in Proposition \ref{lar}
%and conclude that the ideal of the conductor needs at least $11$
%generators.

If the generic $\Gamma'$ is singular at $P$ then it has multiplicity $2$ there.
Then $\overline{D}$ is the blowing-up of $D$ at $P$, denote by $E$
the exceptional divisor. The strict transform $\overline{\Gamma}$ is
base-point-free, because $\overline{\Gamma}$ is semi-ample and
$\overline{\Gamma}^2=8$.

\begin{lem}\label{lemmL} With the assumptions as above the morphism $\rho_D$ does not contract curves.
\end{lem}
\begin{proof} Suppose
that the curve $C\subset D$ is contracted by $\varphi_{|H|}$, then we have
$P\in C$. Thus the image of $C$ is contained in the $3$-dimensional
component of the image of the second exceptional divisor of the Hironaka
model of $(X,H)$ obtained by blowing-up a point then a line in the
exceptional divisor. Now,
since $\overline{\Gamma}|_E$ has degree $2$ the above components
maps to a quadric $Q\subset\mathbb{P}^5$ or to a $3$-dimensional
linear subspace. If the image is linear we can apply the Lemma
\ref{lemm4}. Let us assume that it is a quadric. Then the quadric
$Q$ is contained in a hyperplane $M$. 

In order to end the proof of the Lemma
it is enough to show that $Q$ is
a proper component of $X'\cap M$ cf.~the proof of Lemma \ref{lemm4}. Suppose the contrary, then each
curve $C$ contracted by $\varphi_{|H|}$ is an element of $|\Gamma'|$. Now,
the linear system $|\overline{\Gamma}||_E$ is $2$-dimensional and
$\overline{C}\cdot \overline{\Gamma}=0$. If $\overline{C}\cdot E=m$,
then $$m\overline{\Gamma}\cdot
E=\overline{\Gamma}(\overline{C}+mE)\geq 8$$ so $m\geq 4$. Thus
$|\overline{\Gamma}||_ E$ is $0$-dimensional, a contradiction.
\end{proof}

We deduce that the conductor of the normalization of $X'_D$ defines a
locally CM subscheme $C\subset \mathbb{P}^3$ we also compute as in Section \ref{123} that $\deg (C)=2$.
Now, either $C$ is reduced then it is an aCM
 plane curve or a double
line with Hartshorne--Rao module described in \cite[Ex.~1.5.10]{M}.
We find as in Section \ref{op1}
$$h^1(\mathcal{I}_C(3))=h^1(K_{\overline{D}}-\overline{\Gamma})=6$$
and $h^1(\mathcal{I}_C(n))=0$ for $n>3$ a contradiction.

\subsubsection{Suppose now that $d=7$.} Let us again construct a Hironaka model of $(D,\Gamma')$.
We assume that the support of the base locus is one point $\supp B=\{P\}$. Note that the other cases are simpler and we obtain the same results. 
We consider two possibilities either $\Gamma'$ is smooth or singular at $P$.

If $\Gamma'$ is smooth at $P$
a Hironaka model $\overline{D}$ is obtained by five
successive blow-ups. With the notation as above, the possible
singularities on $Y_D$ are Du Val singularities, rational triple
points, cyclic singularities of type $$\frac{1}{3}(1,1), \ \ \
\frac{1}{4}(1,1),\ \ \  \frac{1}{5}(1,1),$$ and singularities whose
minimal resolutions have exceptional curves with the following
configurations:
\begin{itemize}
\item[]\setlength{\unitlength}{0.4mm}
$$\begin{picture}(70,25)(0,0)
 \put(10,20){\circle{3}}
 \put(20,20){\circle{3}}
 \put(30,20){\circle{3}}
 \put(40,20){\circle{3}}

 \put(7,23){\makebox{\tiny{-4}}}
 \put(17,23){\makebox{\tiny{-2}}}
 \put(27,23){\makebox{\tiny{-2}}}
 \put(37,23){\makebox{\tiny{-2}}}

 \put(11.5,20){\line(1,0){7}}
 \put(21.5,20){\line(1,0){7}}
 \put(31.5,20){\line(1,0){7}}

 \end{picture}
 \begin{picture}(70,)(0,0)
 \put(10,20){\circle{3}}
 \put(20,20){\circle{3}}
 \put(30,20){\circle{3}}
 \put(40,20){\circle{3}}

 \put(7,23){\makebox{\tiny{-3}}}
 \put(17,23){\makebox{\tiny{-2}}}
 \put(27,23){\makebox{\tiny{-2}}}
 \put(37,23){\makebox{\tiny{-3}}}

 \put(11.5,20){\line(1,0){7}}
 \put(21.5,20){\line(1,0){7}}
 \put(31.5,20){\line(1,0){7}}

 \end{picture}
  \begin{picture}(70,)(0,0)
 \put(10,20){\circle{3}}
 \put(20,20){\circle{3}}
 \put(30,20){\circle{3}}

 \put(7,23){\makebox{\tiny{-3}}}
 \put(17,23){\makebox{\tiny{-3}}}
 \put(27,23){\makebox{\tiny{-2}}}

 \put(11.5,20){\line(1,0){7}}
 \put(21.5,20){\line(1,0){7}}

 \end{picture}
 \begin{picture}(70,)(0,0)
 \put(10,20){\circle{3}}
 \put(20,20){\circle{3}}

 \put(7,23){\makebox{\tiny{-4}}}
 \put(17,23){\makebox{\tiny{-2}}}

 \put(11.5,20){\line(1,0){7}}

 \end{picture}
\begin{picture}(70,)(0,0)
 \put(10,20){\circle{3}}
 \put(20,20){\circle{3}}

 \put(7,23){\makebox{\tiny{-3}}}
 \put(17,23){\makebox{\tiny{-3}}}
 \put(11.5,20){\line(1,0){7}}
\end{picture}$$
\end{itemize}
In the figure ``o'' denotes a nonsingular rational curve with
self-intersection equal to the number above it. In each case the
fundamental cycle is equal to the reduced curve and the arithmetic
genus is $0$. Thus by \cite[Thm.~3.5]{A} the singularities on $Y_D$
are rational.

If $\Gamma'$ is singular at $P$ then it has multiplicity $2$ there.
The strict transform of $\Gamma'$ on the blow-up of $D$ at $P$ has
self-intersection $8$, thus it has a base point on the exceptional
divisor. Blowing-up this point we obtain a Hironaka model
$(\overline{D},\overline{\Gamma})$ of $(D,\Gamma')$.

We claim that the morphism given by $|\overline{\Gamma}|$ does not
contract curves. Since $|\overline{\Gamma}|$ maps the exceptional
divisors on $\overline{D}$ into two lines, we can argue as in Lemma
\ref{lemmL}. The claim follows.

In any cases we infer that the subscheme $C$ given by the conductor is locally CM
moreover we compute from  \cite[Thm.~3.5]{KLU} that $\deg (C)<0$ (cf.~Section \ref{123}). This is a contradiction.
\begin{rem} It was observed by O'Grady that the above case $d=7$ and
$\dim B=0$ can be dealt with by comparing the geometric genus of a
generic $$\Gamma'=H_2\cap H_3\cap H_4$$ and the corresponding
birational plane septic curve $$C=L_1\cap L_2\cap L_3\cap X';$$ on
one hand $p_g(\Gamma')\geq 18$ on the other hand $p_g(C)\leq 15$.
\end{rem}

\section{Proof of Theorem \ref{main} in the cases $d\geq9$}\label{rem}

The aim of this section is to give the prove of Theorem \ref{main} in the case 
when the degree $d$ of $X'\subset \PP^5$ is in the range $9\leq d\leq 12$.
In those cases we have to show that if $\varphi_{|H|}$ is birational then the base locus of $|H|$ is zero dimensional.
%Moreover, in the subsections \ref{sec deg11}, \ref{10}, \ref{s9} we study carefully the curve $C\subset\PP^3$ being a generic linear section of the
%singular scheme $\mathcal{C}\subset X'\subset \PP^5$.
%Let us use the notation from section \ref{preliminary}.
\begin{prop} Suppose that $\varphi_{|H|}$ is birational and the degree $d$ of $X'\subset \PP^5$ is in the range $9\leq d\leq 12$. Then the base locus of $|H|$ is $0$ dimensional.
\end{prop}
\begin{proof} It follows from Lemma \ref{4}
 that if $d=11$ then $m=0$. We infer that $\varphi_{|H|}$ has
$0$-dimensional base locus. Moreover, $\supp B$ is exactly one point
$P$ such that $\mult_P(H_1\cdot\Gamma)=1$.

Suppose that $d=10$ that $m\neq 0$ in Lemma \ref{4}. It follows
that $m=1$ and $$\sum_{p\in \supp B}\mult_p (H_1\cdot\Gamma)=0.$$
But this is impossible since the intersection of ample divisors
defining $\Gamma+ \Sigma$ is connected, thus $$\supp\Gamma\cap
\supp\Sigma \neq\varnothing.$$ It follows from Lemma \ref{4}
that $B$ has length $2$.

Let us treat the remaining case $d=9$.
 If $H_1,H_2,H_3\in |H|$ are generic we write $$[H_1\cap H_2\cap
H_3]=\Gamma'+\Sigma'$$ where $\supp\Sigma'\subset \supp B$ and
$$\dim(\supp \Gamma'\cap\supp \Sigma')=0.$$ Denote as above by $D$ the
intersection of $H_1$ and $H_2$ and by $X'_D\subset\mathbb{P}^3$ the
corresponding linear section of $X'\subset \PP^5$. Recall that the surface $D$ is
reduced, irreducible, and normal. Moreover, $D$ locally
Cohen--Macaulay and $\omega_D=2H|_D$.

 Suppose that the base locus has dimension $1$,
i.e.~$\Sigma'\neq 0$. Then by \cite[Prop.~5.4]{O} we deduce that the
cycle $\Sigma'$ is a reduced irreducible local complete intersection
curve and is the scheme-theoretic base locus of $|H||_D$. By
\cite[Thm.~2.1]{DH} the surface $D$ has only isolated singularities
at singular points of $\Sigma'$.
  From the proof of \cite[Prop.~5.4]{O} the curves $$\Gamma'\ \ \ \ \text{and}
\ \ \ \ \Sigma'$$ intersect transversally in one point that varies on $B^1$
when $H_3$ changes (see \cite[(5.3.16)]{O}). It follows that
$\Gamma'$ and $\Sigma'$ are Cartier divisors on $D$. Since the morphism given by the linear system $|\Gamma'|$ is equal
to $\varphi_{|H|}|_D$ (see Remark \ref{poi}) we infer
using Lemma \ref{4} that the linear system $|\Gamma'|$ has no base
points.

\begin{claim}\label{lemm11} With the assumptions above the morphism induced by $\varphi_{|H|}$ on $D$ does not contract
curves.
\end{claim}
\begin{proof} Suppose that the curve $C$ is contracted. Assume moreover, that
$C$ intersects $\Sigma'$ in a smooth point $p$ of $\Sigma'$ so in
particular $\mult_p((\Sigma'+ \Gamma')\cdot C)=\mult_p(H_1\cdot C)$. From Theorem \ref{T}
$cl(C)=mh^3/6$ we have $$2\leq H_1\cdot C=(\Gamma'+\Sigma')\cdot C.$$
Since, $\Gamma'$ and $\Sigma'$ are Cartier divisors we infer $C\cdot
\Sigma'\geq 2$. But $\Gamma'-C$ is effective, this is a
contradiction with $\Sigma'\cdot \Gamma'= 1$.

Let us consider the remaining case where there is a point $P_0\in
\Sigma'$ such that for a generic choice of $D$ we have that $P_0\in
C$, where $C$ is a contracted curve. Let $$\pi:\tilde{X}\rightarrow
X$$ be the blow up of $\Sigma'\subset X$. Note that $\tilde{X}$ can have
singularities at points in the pre-image by $\pi$ of singular points
of $\Sigma'$. Now each fiber of the exceptional divisor
$E\rightarrow \Sigma'$ map by $|\pi^{\ast}(H)-E|$ to a plane in
$\mathbb{P}^5$. It follows that the image of curves contracted by
$\varphi_{|H|}$ is contained in a $2$-dimensional linear subspace of
$\mathbb{P}^5$ being the image of the fiber of $E$ that maps to $P_0$.
We conclude by Lemma \ref{lemm4}.\end{proof}

From \cite[Prop.~2.3]{R}, we infer $$5\Gamma'=2\Gamma'+2\Sigma' + C_1$$
where $C_1\subset D$ is the Cartier divisor defined by the
conductor. Finally, using \cite[Thm.~3.5]{KLU}, we compute that the
degree of the conductor subscheme $C\subset X_D'\subset\mathbb{P}^3$
is $$\frac{1}{2}\Gamma'(3\Gamma'-2\Sigma')=\frac{25}{2}$$ and obtain
a contradiction.
\end{proof}
We have completed the proof of Theorem \ref{main}.

\section{The singular locus of $X'\subset \PP^5$ in the cases $d\geq 9$}\label{last}

Our aim is to describe $X'\subset \PP^5$ in the case when $\varphi_{|H|}$ is birational and $d\geq 9$.
Let us first treat the cases $d=9,10$. In the case $d=11$ the map $\varphi_{|H|}$ does not contract curves; we consider it separately in Section \ref{sec deg11}.

From the previous section we can assume that the base locus of $|H|$ is $0$ dimensional.
By Lemma \ref{lemmO} we can assume that $D$ the intersection of two divisors from $|H|$ is smooth.
In order to construct the normalization of $\varphi_{|H|}(D)=X'_D\subset \PP^3$ we need the following:

\begin{prop}\label{contr-D} If $\varphi_{|H|}$ is birational and $d=9$ or $10$ then the morphism induced by $\varphi_{|H|}$ on $D$ does not contract
curves.
\end{prop}
\begin{proof} 
Assume that $d=10$.
It is enough to prove that $\beta(\rho(E))$ is a sum of linear
spaces of $\dim\leq 3$ i.e.~ that the assumptions of Lemma \ref{lemm4} holds. 

Let us construct a Hironaka models.
 Suppose that $\supp B$ is one point since the case when $\supp B$ are two points is clear.
After the first blow up $\gamma$ the system $$|\gamma^{\ast}(H)-E_1|$$ restricted to ${E_1}$ is the system of
hyperplanes passing through one point $Q$. The pair $(\overline{X},H_2)$,
where $H_2$ is the proper transform of $H$, obtained by the
blowing-up of $Q$, is a Hironaka model of $(X,H)$. Thus the exceptional divisors maps to linear spaces of dimension $\leq 3$.

%From Lemma \ref{lemm2} the restriction of the system $|\gamma^{\ast}(H)-E_1|$ to ${E_1}$ is a linear
%subsystem of $|\mathcal{O}_{E_1}(1)|$. Thus it is a linear system of
%hyperplanes passing through a linear subspace $\Lambda\subset E_1$.
%We claim that the dimension of $\Lambda$ is $0$. Indeed, suppose it
%is larger. Then $\Gamma$ has to be singular. Since $\Gamma$ is tangent to $H_1$ it follows that
%$$\mult_p(H_1\cdot\Gamma)\geq 4.$$ 
%Now, since $$\mult_p(H_1\cdot\Gamma)=2$$ the blowing-up of $Q$
%separates the proper transforms of $\Gamma$ and $H_1$.
%The claim follows.

Finally consider the case $d=9$.
If $\supp B $ are two or three points, we argue as in the cases of degrees
$10$. If $\supp B$ is one point $P$ then each contracted
curve contain $P$, 
the proof is analogous as for $d=9$. Note that a Hironaka model is constructed in Section \ref{afa}. 
%we conclude as in the proof of Lemma \ref{lemm4}.
This finish the proof.
\end{proof}

We shall now treat the cases $d=9$ and $ 10$ separately.
The case $d=11$ is discussed in the last section.

\subsection{Degree 9}\label{s9}

Let us define the surface $D$ as above as the intersection of two divisors from $|H|$. From Lemma \ref{lemmO} the
surface $D$ is smooth and from Proposition \ref{contr-D} the morphism $\varphi_{|H|}$ do not contract divisor on $D$.
The aim of this section is to find the possible cohomology table of the singular locus of $X_D'\subset \PP^3$ and give evidence that such curve cannot exist.

\subsubsection{A Hironaka model of $X'_D\subset \PP^3$}\label{afa} Let us assume that $\Supp B$ is one point $P$. In the other cases we obtain the same cohomology table for the ideal sheaf of $C\subset \PP^3$.

%Thus the linear system $|\Gamma|$ on $D$ can have only some base points.
Let $D'$ be the surface obtained from $D$ by blowing-up $P$.
 Denote by $\Gamma''$ the strict transform of $\Gamma'$ on $D'$.
We show as in the proof of Lemma \ref{lemmO} that the generic
element of $|\Gamma'|$ is smooth at $P$. Thus $$(\Gamma'')^2=11$$ and
$\Gamma''$ has exactly one base point $P'$ on the exceptional
divisor $E'$ on $D'$, moreover $\Gamma''\cdot E'=1$. Blowing-up $P'$
we obtain a surface $D''$ with exceptional divisor $E''$ such that
$\Gamma'''$ (resp.~$E''$) is the strict transform of $\Gamma''$
(resp.~$E'$). We have $$(\Gamma''')^2=10,$$ the linear system
$|\Gamma'''|$ has exactly one base point $P''\in E''$, and
$\Gamma'''\cdot E''=1$. The strict transform $\overline{\Gamma}$ of
$\Gamma'$ on $\overline{D}$, the blowing-up of $D''$ at $P''$ gives
a base-point-free linear system; denote by
$$\pi_D:\overline{D}\rightarrow D$$ the composed morphism. The
morphism $\rho$ defined by $|n\overline{\Gamma}|$ for $n$ large
enough has normal image $Y_D$ and contracts only the strict
transforms of $E'$ and $E''$. It follows from the Stein
factorization theorem that $Y_D$ is the normalization of the chosen
codimension $2$ linear section $X_D'$ of $X'\subset \mathbb{P}^5$.

We have two possibilities: either $$P''\in E''- E' \ \ \ \  \text{or} \ \ \ P''\in E'.$$
First observe that the possibility $P''\in E'$
cannot happen from the fact that the base scheme is a
curvilinear scheme of length $3$ (see Remark \ref{lp}).

In the case $P''\in E''- E'$ we infer from the Artin contraction theorem that
$Y_D$ has exactly one singular point being a Du Val singularity of
type $A_2$. In particular, $Y_D$ is locally Cohen--Macaulay, and
$\omega_{Y_D}$ is locally free. 

\subsubsection{The singular locus of $X_D'\subset \PP^3$}\label{123}
The ideal of the conductor of the normalisation defines a
subscheme $C\subset X'_D$ that is of pure dimension $1$ and locally
Cohen--Macaulay (from the proof of \cite[Thm.~3.1]{Ro}).
Finally, using \cite[Thm.~3.5]{KLU} and  \cite[Prop.~2.3]{R}, we compute the degree
$$2deg(C)=L(5L-K_{Y_D})=45-\rho^{\ast}(K_{Y_D})\cdot \overline{\Gamma}=45-\overline{\Gamma}
\cdot (K_{\overline{D}}+aE'+bE'')=$$ $$=45-\overline{\Gamma}\cdot
K_{\overline{D}}=45-(\pi_D^{\ast}(2\Gamma')+E'+2E''+3F)(\pi_D^{\ast}(\Gamma')-E'-2E''-3F)=18.$$
 We compute also
$$\chi(K_{Y_D}+nL)=\frac{9}{2}n^2+\frac{27}{2}n+12$$ and $p_a(C)=1$.

%Next, if $P''\in E'$ then $E'$ and $E''$ are disjoint $-3$ and $-2$
%curves on $\overline{D}$. We find that $Y_D$ has exactly two
%singular points: a quotient singularity of type $\frac{1}{3}(1,1)$
%(see \cite[Rem.~4.9]{KM})
% and a Du Val singularity of type $A_2$. It follows from \cite[Prop.~5.15]{KM} that $Y_D$
%has rational singularities (thus is locally Cohen--Macaulay) and is
%$\mathbb{Q}$-factorial. We infer that the conductor of the
%normalization $\beta_D: Y_D\rightarrow X'_D$ defines a subscheme
%$C\subset X'_D$ that is of pure dimension $1$ and locally
%Cohen--Macaulay.
\begin{prop}\label{lar}  We have $h^0(\mathcal{I}_C(3))=1$,
$h^0(\mathcal{I}_C(4))=4$, $h^0(\mathcal{I}_C(5))=11$, and
$h^0(\mathcal{I}_C(6))=30$. Moreover, $C$ can be $3\times 6$
linked to $B\subset
\mathbb{P}^3$ a degree $9$ curve with $p_a(B)=1$.
\end{prop}
\begin{proof}
As before we deduce from (\ref{2}) that
$$h^1(\mathcal{I}_C(n+5))=h^1((\beta_D)_{\ast}(\omega_{Y_D})(n)).$$
Define $L:=\beta_D^{\ast}(\mathcal{O}_{X'_D}(1))$. Then since $Y_D$
has rational singularities, we infer from \cite[Cor.~6.11]{EV} that
$$h^i(K_{Y_D}+nL)=h^i(K_{\overline{D}}+n\overline{\Gamma}).$$  Now, let us compute
$$h^0(K_{\overline{D}}+n\overline{\Gamma})-h^0(\mathcal{O}_{\mathbb{P}^3}(-4+n))=h^0(\mathcal{I}_C(n+5)).$$
Denote by $F$ the exceptional divisor of the last blow-up. Then
$$K_{\overline{D}}= \pi_D^{\ast}(K_D)+E'+2E''+3F$$ and
$\pi_D^{\ast}(\Gamma')=\overline{\Gamma}+E'+2E''+3F.$

We claim that
$$h^0((n+2)\overline{\Gamma}+3E'+6E''+9F)=h^0(K_{\overline{D}}+n\overline{\Gamma})=h^0(K_{D}+n\Gamma')$$ for
$n=-2,-1,0,1$. From \cite[Lem.~4.3.16]{L} we have
$$ h^0((n+2)\overline{\Gamma})=h^0((n+2)\Gamma')=h^0((n+2)(\overline{\Gamma}+E'+2E''+3F)).$$
The claim follows. 
%using the long exact sequence as in Section\ref{sec deg11}.

 Finally we have the following:
 \begin{center}
            %\begin{table}
            \renewcommand*{\arraystretch}{1}
\begin{table}[htp]
 ×

\begin{tabular}{|c|c|c|c|c|}

\hline
            $n$& $h^0(K_{Y_D}+nL)$& $h^1(K_{Y_D}+nL)$& $h^2(K_{Y_D}+nL)$& $\chi(K_{Y_D}+nL)$                       \\ \hline
                                                                   \hline

            $0$&11&$0 $&$1$& $12$           \\ \hline
            $-1$&$4$ &$5$&$4$& $3$ \\ \hline
            $-2$& 1&$y$&$2+y$&$3$                    \\ \hline
            $-3$&0&$z$&$12+z$&$12$\\ \hline
            $-4$&$0$&$t$&$30+t$&$30$\\ \hline
            \end{tabular}
 \end{table}
           %\end{table}
            \end{center}
We see that $9\geq y\geq 8$, $z\geq 8$, and $t\geq5$. This finish the proof. 
\end{proof}

Let us study the the curve $B\subset
\mathbb{P}^3$ in order to find evidence to the fact that $C$ cannot exist.
\begin{center}
            %\begin{table}
            \renewcommand*{\arraystretch}{1}
\begin{table}[htp]
 ×

\begin{tabular}{c|ccc}

            $n$& $h^0(\mathcal{I}_B(n))$& $h^1(\mathcal{I}_B(n))$& $h_B(n)$                       \\ \hline

            $4$&$t-1$ &$t$&$0$ \\
            $3$&$z-7$ &$z$&$0$ \\
            $2$& $y-8$&$y$&$1$                    \\
            $1$&0&$5$&$7$

             \end{tabular}
 \end{table}
           %\end{table}
            \end{center}

First if $y=9$ then $h_C(5)=1$, it follows from \cite[Thm.~1.1]{Sl}
that $h_C(4)=h_C(3)=1$ (see the $1$-property \cite{S3}), thus $z=11$
and $t=11$.
 We can use \cite[Cor.4.4]{S2} to show that $B$ is
not minimal. But $B$ can be bilinked down with height $-1$ on the
quadric to the minimal curve $B_0$ of degree $7$. We compute form
Equation (\ref{leq3}) that $h^0(\mathcal{I}_{B_0}(2))=1$ and
$h^0(\mathcal{I}_{B_0}(3))=4$. The Betti table of the minimal
resolution of $M_B$ is as follows.
\begin{center}
            %\begin{table}
            \renewcommand*{\arraystretch}{1}
\begin{table}[htp]
 ×

\begin{tabular}{c|cccc}

            $j\setminus i$& $0$& $1$& $2$  &\dots                     \\ \hline

            $1$&5& $11$ &$5+k$&\dots  \\
            $2$& 0&$k$&$d$ &\dots                    \\
            \dots &$0$ &0 &\dots&\dots

            \end{tabular}
 \end{table}
           %\end{table}
            \end{center}
With the notation of section \ref{liaison} we infer that $q(2)=4+k$ thus $q(3)=1$ since $h=0$.

 Assume that
$y=8$, the Betti table of the minimal resolution of $M_B$ is as
follows.
\begin{center}
            %\begin{table}
            \renewcommand*{\arraystretch}{1}
\begin{table}[htp]
 ×

\begin{tabular}{c|cccc}

            $j\setminus i$& $0$& $1$& $2$  &\dots                     \\ \hline

            $1$&5& $12$ &$z-2+k$&\dots  \\
            $2$& 0&$k$&$d$ &\dots                    \\
            \dots &$0$ &0 &\dots&\dots

            \end{tabular}
 \end{table}
           %\end{table}
            \end{center}
If $B$ is minimal in its biliaison class then $q(2)=5+k$ thus
$q(3)=1$ since $h=1$. We can find a bound of the invariant $a_1$
from \cite[p.~77]{MP}. This gives an evidence for the Conjecture \ref{conj} in the case $d=9$
since $a_1$ is different than expected.

\begin{lem} The invariant $a_1>3$. \end{lem}
\begin{proof} Suppose the contrary i.e. $a_1=3$ then from \cite[Prop.~5.10IV]{MP} $B$ is $3\times 4$ linked to a minimal curve $C_0$
from the class of $C$. Then $$\deg C_0=3, \ \ \ p_a(C_0)=-8,\ \ \ 
h^0(\mathcal{I}_{C_0}(2))=0, \ \ \text{and}\ \ \ h^0(\mathcal{I}_{C_0}(3))=2.$$ Since
$C_0$ is not extremal we find $e(C_0)=-2$ thus $C_0$ has a
quasi-primitive structure supported on a sum of lines
\cite[Rem.~3.5]{Sl} and non reduced \cite[Ex.~2.11]{Sl}. If $C_0$ is
supported on two lines then we obtain a contradiction with
$h^0(\mathcal{I}_{C_0}(3))=2$ from the proof of \cite[Prop.~3.3]{N1}
(and by \cite[Prop.~3.2]{N1} since $C_0$ is not extremal). If it is
supported on one line then the possible number of cubic generators
of $\mathcal{I}_{C_0}$ are computed in \cite[Rem.~2.4,
Prop.~2.1]{N1} this is a contradiction.
\end{proof}
If $B$ is not minimal it can be bilinked down (on a cubic) to a
minimal curve $B_0$. From (\ref{leq3}) we deduce that $q(2)\geq 6+k$
thus $h=0$, $q(2)=6+k$, and $a_1>3$.
\begin{prob} A curve $C\subset \PP^3$ of degree $9$ with the invariants described above does not exist. 
\end{prob}
%---------

\subsection{Degree 10}\label{10} Let us study the curve $C\subset \PP^3$ in the case $d=10$.
 Let us construct the normalization $Y_D$ when $\Supp(D)$ is one point since the case when $\Supp(D)$ is two points is simpler and give the same cohomology table for $C\subset \PP^3$. Recall that a Hironaka model is constructed in the proof of Proposition \ref{contr-D}. So we need only the following:
\begin{lem} The morphism $\rho_D$ is $1:1$ on $\overline{D}-E $, where $E$ is the strict transform of the exceptional
curve of the first blow-up on $\overline{D}$. Moreover, $\rho_D$
contracts $E$ to a Du Val singularity of type $A_1$ on $Y_D$.
\end{lem}
\begin{proof}
 The system $|H'_D|$ being the restriction of $|H|$ to $D$ does not contract
curves. Now, $E$ is a smooth rational curve with self intersection
$-2$, that is contracted by $\rho_D$ to a normal singularity. This
singularity must be an ordinary double point.
\end{proof}

 We deduce that the normalization $Y_D$ is locally Cohen--Macaulay and $\omega_{Y_D}$ is
locally free. From Theorem \ref{Za} the conductor $C\subset X'_D$
has pure dimension $1$. Since $K_{Y_D}$ is a Cartier divisor, and
$$C\in |6L-K_{Y_D}|$$ we deduce that $C$ is locally Cohen--Macaulay. 
Denote $L:=\beta
^{\ast}(\mathcal{O}_{X'_D}(1))$. Then we
infer
$$\beta_D^{\ast}(\mathcal{O}_{X'_D}(6))=\beta_D^{\ast}(K_{X'_D})=K_{Y_D}+C_1,$$ where $C_1\subset Y_D$
is the Cartier divisor defined by the conductor.

Since $Y_D$ is
Cohen--Macaulay, 
%the normalization is locally flat of codimension $1$ 
we compute as in the case $d=9$ that $2\deg(C)=\deg(C_1)$.
Now, $Y_D$ has rational singularities and it is $\mathbb{Q}$-factorial,
thus we can compute as follows $\deg(C_1)=L(6L-K_{Y_D})$. Since $H_D\cdot K_{\overline{D}}=26$ we have
$$2\deg(C)=6H_D^2-H_D\rho_D^{\ast}(K_{Y_D})=34.$$  Moreover,
 we also obtain $p_a(C)=30$ and arguing as before the following:

\begin{prop} The curve $C$ is contained in one quartic such that  $h^0(\mathcal{I}_C(4))=1$,  $h^0(\mathcal{I}_C(5))=4$,
$h^0(\mathcal{I}_C(6))=11$, $h^0(\mathcal{I}_C(7))=30$ and $C$ is generated in degree $7$. Moreover, $C$ can be $4\times 7$ linked to a curve of degree $11$.
\end{prop}

\subsection{Degree 11} \label{sec deg11}  In this case the methods of proof of Theorem \ref{main} does not work.
%n order to study the non-normal locus of $X_D'$ we shall construct the normalization $Y_D\to X'_D$; in each degree separately.

 % Moreover, t that if $y=2$, $z=t=\dts=0$ then
In fact L. Gruson pointed out that 
a general curve $C\subset \PP^3$ with the required numerical type exists and can be obtained by starting
with a general curve $C_0$ of degree $9$ and genus $6$ and taking the general linked curve by
surfaces of degrees $(5, 7)$. The ideal of $C$ has the presentation
                         $$0\rightarrow \mathcal{O}_{\mathbb{P}^3}^{\oplus 4} (-8)\rightarrow \mathcal{O}_{\mathbb{P}^3}(-7)\rightarrow \mathcal{O}_{\mathbb{P}^3}(-5)\oplus E(-4)$$
  where $E$ is a well-defined ``Schwarzenberger'' bundle of rank 3. Note that $p_a(C)=74$ in this case.

  Moreover, he observed that $\varphi_{|H|}$ does not contract curve on $X$ when $d=11$. Thus
  we do not need to consider the codimension $2$ section $X_D'$ in order to construct the normalization.

\begin{lem}
 If $d=11$ the map $\rho$ induced by $\varphi_{|H|}$ from $Z$ the blown-up of $X$ in the point of indeterminacy $P$, to $X'\subset \PP^5$ is finite.
\end{lem}
\begin{proof}
 Suppose that a curve $C$ is contracted by $\varphi_{|H|}$. Since $H$ is ample we have that $C.H>0$ and from Thm.~1.1(3) we have 
$C.H$ is even. Since $C$ is contracted by $|H|$, it follows that $$(C.H)_P>1.$$ Now, the strict transform $\overline{C}$ of $C$ on $Z$ cuts the exceptional 
divisor $E$ with multiplicity greater then one. Moreover the system $$|\overline{H}|=|\pi^*(H)-E|$$ is base-point-free thus on $E$ restricts to the complete linear system of $\mathcal{O}_E(1)$. Note that $E$ is isomorphic to 
$\mathbb{P}^3$ and a generic element of $|H|$ is smooth.        
We obtain a contradiction if $\overline{C}$ cuts $E$ in two different points of $E$.

Assume that $\overline{C}$ cuts $E$ in one point with multiplicity $>1$.
Localy the map $\varphi_{|\overline{H}|}$ can be seen as follows: 
$$\mathbb{C}^4\supset U\rightarrow \mathbb{C}^5$$ given by holomorphic functions $f_1,\dots,f_5$ vanishing on $\overline{C}\cap U$ and maping $P=0$ to $0$ ($U$ is an open subset).
Consider the functions $f_i|_E=g_i$ for $i=1,\dots,5$. Since $\overline{C}$ cuts $E$ with multiplicity $>1$ in $P$, we obtain that the functions $g_i$ defines on 
$E$ a non reduced point with support $P$.
It follows that the rank of the differential of $(g_1,\dots,g_5)$ cannot have maximal rank at $0$ this is a contradiction since it has to be an isomorphism. \end{proof}
We infer as before the following:
\begin{lem} There is a unique quintic hypersurface containing the singular locus of $X'\subset \PP^5$.
\end{lem}
Let us describe the image $X'\subset \PP^5$ more precisely.
We can apply the Beilinson monad for $E_1^{pq} = \Omega^p (p) \otimes H^q(\oo_Z (2-p))$.
First we compute that $$h^i(\mathcal{O}_{X'}(n))=h^i(\mathcal{O}_Z(n\overline{H}))=:h^i(\mathcal{O}_Z(n))$$ for $i=2,3,4$ and $2\geq n\geq-3$ from the Leray spectral sequence. 
%for the blow-up and the long cohomology exact sequence
 %deduced from $0\rightarrow \mathcal{O}_Z(-E)\rightarrow\mathcal{O}_Z\rightarrow \mathcal{O}_E\rightarrow 0)$.
Next we have an exact sequence \[H^0(\mathcal{O}_Z(\overline{H}))\xrightarrow{\alpha} H^0(\mathcal{O}_Z(\overline{H}+E))\rightarrow H^0(\mathcal{O}_E)\rightarrow H^1(\mathcal{O}_Z(\overline{H}))\rightarrow H^1(\mathcal{O}_Z(\overline{H}+E))\rightarrow0.\]

Now the map $\alpha$ is a surjection since $P$ is a simple base point.
We deduce that $$H^0(\mathcal{O}_{E})= H^1(\mathcal{O}_Z(\overline{H}))$$ and $h^1(\mathcal{O}_Z(n))=h^1(\mathcal{O}_X(n))$ for $0\geq n\geq-3$. We identify $\mathbb{P}^3$ with the image of the exceptional divisor in $\mathbb{P}^5$.
Next using the fact that $H^0(\mathcal{O}_X(2H))=Sym^2(H^0(\mathcal{O}_X(H)))$ we deduce that $P$ is a simple base point of $|2H|$ thus 
$$H^0(\mathcal{O}_{E})\oplus H^0(\mathcal{O}_{E}(1))= H^1(\mathcal{O}_Z(2\overline{H}))$$ and $h^0(\mathcal{O}_Z(n))=h^0(\mathcal{O}_X(n))$ for  $2\geq n\geq-3$.

\begin{center}\begin{tabular}{cccccc}
$ H^4(\mathcal{O}_X(-3))$& $ H^4(\mathcal{O}_X(-2))$&$ H^4(\mathcal{O}_X(-1))$&
$\mathbb{C} $&$0$&$0$\\
$0$&$0$&$0$&$0$&$0$&$0$\\
$0$&$0$&$0$&$\mathbb{C} $&$0$&$0$\\
$0$&$0$&$0$&$0$&$H^0(\mathcal{O}_{\PP^3})$&$H^0(\mathcal{O}_{\PP^3})\oplus H^0( \oo_{\PP^3}(1))$  \\
$0$&$0$&$0$&$\mathbb{C} $&$H^0(\mathcal{O}_X(1))$&$ H^0(\mathcal{O}_X(2))$
\end{tabular}
\end{center}
Now arguing as in \cite{CaS} (see \cite[\S 4]{K}) we infer that the Beilinson monad applied to $\rho_{\ast}\oo_Z(2)$ associated to the above table is cohomologous to the following:
\begin{center}
\begin{tabular}{cccccc}
$\Omega_{\PP^5}^5(5)\otimes A\oplus \oo_{\PP^5}(-4)$ &$0$&$\ \ \ \ \ \ \ 0\ \ \ \ \ \ \ \ $&$\ \ \ \ \ \ \  0\ \ \ \ \ \  $&$\ \ \ \ \ \ 0\ \ \ \ \ \ $&$0$\\
$0$&$0$&$0$&$0$&$0$&$0$\\
$0$&$0$&$0$&$\Omega_{\PP^5}^2(2)$&$0$&$0$\\
$0$&$0$&$0$&$0$&$\Omega^1_{\PP^5}(1)$&$5\oo_{\PP^5}$  \\
$0$&$0$&$0$&$0$&$0$&$\mathcal{O}_{\mathbb{P}^5}(2)$
\end{tabular}
 \end{center}
So the monad have the following shape:
$$0\to \oo_{\PP^5}(-4)\oplus 10\oo_{\PP^5}(-1)\to \Omega^2_{\PP^5}(2)\oplus \Omega^1_{\PP^5}(1)\oplus \oo_{\PP^5}(2) \to5\oo_{\PP^5}\to 0, $$
such that $\rho_{\ast}(\mathcal{O}_Z(2))$ (with support $X'$) is the cohomology in the middle. 
This will be used in a future work to study the geometry of $X'$.

%Moreover, in the monade the map $$d_1^{-1,1}\colon \Omega_{\mathbb{P}^5}^1(1)\otimes H^0(\mathcal{O}_{\mathbb{P}^3})\rightarrow \mathcal{O}_{\mathbb{P}^5}\otimes( H^0(\mathcal{O}_{\mathbb{P}^3}(1))\oplus H^0( \mathcal{O}_{\mathbb{P}^3}))$$
%is induced by the natural multiplication map $H^0(\mathcal{O}_{\mathbb{P}^3}) \otimes H^0(\mathcal{O}_{\mathbb{P}^5}(1)) \rightarrow H^0(\mathcal{O}_{\mathbb{P}^3}(1))$.
%Remark that $d_1^{-1,1}|_{\mathbb{P}^3}$ is induced from the Euler sequence thus $\ker(d_1^{-1,1})$ should be a line bundle. Computing the Chern classes we obtain
%$$\ker(d_1^{-1,1})=\mathcal{O}_{\mathbb{P}^5}(-1)$$  $$\coker(d_1^{-1,1})=\mathcal{O}_{\mathbb{P}^3}(1) \oplus \mathcal{O}_{\mathbb{P}^5}.$$
%Finally from the Belinson theorem the map $d_2^{0,2}\colon \Omega_{\mathbb{P}^5}^2(2)\to \mathcal{O}_{\mathbb{P}^3}(1) \oplus \mathcal{O}_{\mathbb{P}^5}$ is surjective denote by $T$ its kernel (it is a coherent sheaf). 
%We obtain an exact sequence $$0\to\mathcal{O}_{\mathbb{P}^5}(-4)\oplus 10 \mathcal{O}_{\mathbb{P}^5}(-1)\to T\oplus \mathcal{O}_{\mathbb{P}^5}\oplus\mathcal{O}_{\mathbb{P}^5}(2)\to f_{\ast}(\mathcal{O}_Z(2))\to 0.$$
%\begin{rem}Denote $U=\mathbb{P}^5-\mathbb{P}^3$ the complement of the exceptional divisor. Denote by $T_U$ the vector bundle being the kernel of the surjection $\Omega_U^2(2)\to \mathcal{O}_U$ from the monade.
%Is the quintic in the ideal of the singular locus of $\varphi_{|H|}(X)$ the closure of the determinantal locus of a injection $9\mathcal{O}_U\to T_U$?
%\end{rem}

%---------------
           
\renewcommand{\refname}{References:}

\end{document}